\documentclass{gtart}

\def\ifplaintex{\expandafter\ifx\csname documentclass\endcsname\relax}

\ifplaintex 
\else
\fi

\expandafter\ifx\csname beginpicture\endcsname\relax
\expandafter\ifx\csname documentclass\endcsname\relax
\input pictex \else
\input prepictex \input pictex \input postpictex \fi\fi

\def\gt{{\mathsurround=0pt\it $\cal G\mskip-2mu$eometry \&\ 
$\cal T\!\!$opology}}        

\def\gtp{{\mathsurround=0pt\it $\cal G\mskip-2mu$eometry \&\ 
$\cal T\!\!$opology $\cal P\!$ublications}}  

\def\lognumber#1{\def\thelognumber{#1}}
\def\volumenumber#1{\def\thevolumenumber{#1}}
\def\papernumber#1{\def\thepapernumber{#1}}
\def\volumeyear#1{\def\thevolumeyear{#1}}

\def\pagenumbers#1#2{\def\startpage{#1}\def\finishpage{#2}}
\def\published#1{\def\publishdate{#1}}
\def\proposed#1{\def\theproposer{#1}}
\def\seconded#1{\def\theseconders{#1}}
\def\received#1{\def\receiveddate{#1}}
\def\revised#1{\def\reviseddate{#1}}
\def\accepted#1{\def\accepteddate{#1}}
\def\asciititle#1{\def\theasciititle{#1}}

\def\asciiaddress#1{\def\theasciiaddress{#1}}
\def\asciiemail#1{\def\theasciiemail{#1}}
\long\def\asciiabstract#1{\long\def\theasciiabstract{#1}}

\let\\\par\let\thelognumber\relax
\let\thevolumenumber\relax\let\thepapernumber\relax
\let\thevolumeyear\relax\let\thesamplenumber\relax\let\startpage\relax
\let\finishpage\relax\let\publishdate\relax\let\receiveddate\relax
\let\reviseddate\relax\let\accepteddate\relax\let\theasciititle\relax
\let\theasciiauthors\relax\let\theasciiaddress\relax
\let\theasciiabstract\relax
\let\theasciiemail\relax\let\theshortauthors\relax\let\theshorttitle\relax

\long\def\maketitlep{   

\count0=\startpage

\gt\hfill      
\beginpicture
\setcoordinatesystem units <0.33truein, 0.33truein> point at 2.2 0.9
\setplotsymbol ({$\cal G$})
\plotsymbolspacing=9truept
\circulararc 315 degrees from 0 1 center at 0 0
\setplotsymbol ({$\cal T$})
\circulararc 315 degrees from 1 -1 center at 1 0
\endpicture
\break
{\small\ifx\thesamplenumber\relax 
Volume \else Sample
\fi\thevolumenumber\ (\thevolumeyear)
\startpage--\finishpage\nl
Published: \publishdate}
\vglue 0.5truein plus 0.4fil minus 0.1truein

{\parskip=0pt\leftskip 0pt plus 1fil\def\\{\par\smallskip}{\ifplaintex\large
\else\Large\fi\bf\thetitle}\par\medskip}   

\vglue 0pt plus 0.1fil 

{\parskip=0pt\leftskip 0pt plus 1fil\def\\{\par}{\sc\theauthors}
\par\medskip}

\vglue 0pt plus 0.1fil 

{\small\parskip=0pt\let\newline\\
{\leftskip 0pt plus 1fil\def\\{\par}{\sl\theaddress}\par}
\expandafter\ifx\theemail\relax    
\relax\else\vglue 5pt plus 0.02fil minus 2pt\def\\{\stdspace{\rm 
and}\stdspace} 
\cl{Email:\stdspace\tt\theemail}\fi
\ifx\theurl\relax                  
\relax\else\vglue 5pt plus 0.02fil minus 2pt\def\\{\stdspace{\rm 
and}\stdspace}
\cl{URL:\stdspace\tt\theurl}\fi\par}

\vglue 7pt plus 0.3fil minus 3pt

{\bf Abstract}
\vglue 5pt plus 0.1fil minus 2pt

\theabstract

\vglue 7pt plus 0.3fil minus 3pt

{\bf AMS Classification numbers}\quad Primary:\quad \theprimaryclass

Secondary:\quad \thesecondaryclass

\vglue 5pt plus 0.3fil minus 2pt

{\bf Keywords}\quad \thekeywords

\vglue 10pt plus 0.5fil minus 5pt

{\small  Proposed: \theproposer\hfill Received: \receiveddate\nl
Seconded: \theseconders\hfill 
\ifx\reviseddate\relax                         
Accepted: \accepteddate                        
\else
Revised: \reviseddate                          
\fi}
\eject
}       

\let\maketitlepage\maketitlep
\let\maketitle\maketitlepage

\font\phead=cmsl9 scaled 950
\font=cmsl9 scaled 1050
\font\pnum=cmbx10 scaled 913
\font\lnum=cmbx10 
\font\pfoot=cmsl9 scaled 950
\font=cmsl9 scaled 1050
\ifplaintex
\headline{\vbox to 0pt{\vskip -4.5mm\line{\small\phead\ifnum
\count0=\startpage ISSN 1364-0380 (on line)
1465-3060 (printed) \hfill {\pnum\folio}\else\ifodd\count0\def\\{ }%
\ifx\theshorttitle\relax\thetitle\else\theshorttitle\fi\hfill{\pnum\folio}
\else\def\\{ and }{\pnum\folio}\hfill\ifx\theshortauthors\relax\theauthors
\else\theshortauthors\fi\fi\fi}\vss}}
\footline{\vbox to 0pt{\vglue 0mm\line{\small\pfoot\ifnum\count0=\startpage
\copyright\ \gtp\hfill\else
\gt, Volume \thevolumenumber\ (\thevolumeyear)\hfill\fi}\vss
}}
\else
\makeatletter
\def\@oddhead{{\small\lhead\ifnum\count0=\startpage ISSN 1364-0380 (on line)
1465-3060 (printed) \hfill {\lnum\number\count0}\else\ifodd\count0
\def\\{ }\ifx\theshorttitle\relax \thetitle \else\theshorttitle\fi\hfill
{\lnum\number\count0}\else\def\\{ and }{\lnum\number\count0}
\hfill\ifx\theshortauthors\relax 
\theauthors\else\theshortauthors\fi\fi\fi}}\def\@evenhead{\@oddhead}
\def\@oddfoot{\small\lfoot\ifnum\count0=\startpage\copyright\ \gtp\hfill\else
\gt, Volume \thevolumenumber\ (\thevolumeyear)\hfill\fi}
\def\@evenfoot{\@oddfoot}
\makeatother
\fi

\newwrite\gtoutfile
\long\gdef\makeheadfile{  
{\def\\{, }\def\s{ }
\immediate\openout\gtoutfile head.xxx
\immediate\write\gtoutfile{To: math@arxiv.org}
\immediate\write\gtoutfile{Subject: put or rep NNNNN:pppp}
\immediate\write\gtoutfile{--text follows this line--}
\immediate\write\gtoutfile{Proxy-for: \ifx\theasciiauthors\relax
\theauthors\else\theasciiauthors\fi\s<\ifx\theasciiemail\relax\theemail\else\theasciiemail\fi>}
\immediate\write\gtoutfile{\noexpand\\}
\immediate\write\gtoutfile{Authors: \ifx\theasciiauthors\relax
\theauthors\else\theasciiauthors\fi}
\immediate\write\gtoutfile{Title: \ifx\theasciititle\relax
\thetitle\else\theasciititle\fi}
\immediate\write\gtoutfile{Subj-class: GT or SG or MG etc}
\immediate\write\gtoutfile{MSC-class: \theprimaryclass\ifx\thesecondaryclass\relax\else, \thesecondaryclass\fi}
\immediate\write\gtoutfile{Journal-ref: Geom. Topol. \thevolumenumber
(\thevolumeyear) \startpage-\finishpage}
\immediate\write\gtoutfile{Comments: Published by Geometry and Topology at}
\immediate\write\gtoutfile{\s\s http://www.maths.warwick.ac.uk/gt/GTVol\thevolumenumber/paper\thepapernumber.abs.html}
\immediate\write\gtoutfile{\noexpand\\}
\immediate\write\gtoutfile{}
\ifx\theasciiabstract\relax
\immediate\write\gtoutfile{\theabstract}\else
\immediate\write\gtoutfile{\theasciiabstract}\fi
\immediate\write\gtoutfile{}
\immediate\write\gtoutfile{\noexpand\\}
\immediate\write\gtoutfile{}
\immediate\closeout\gtoutfile}}  

\def\maketitlepage{\maketitlep\makeheadfile}
\let\maketitle\maketitlepage

\def\ifplaintex{\expandafter\ifx\csname documentclass\endcsname\relax}

\ifplaintex 
\else
\fi

\expandafter\ifx\csname beginpicture\endcsname\relax
\expandafter\ifx\csname documentclass\endcsname\relax
\input pictex \else
\input prepictex \input pictex \input postpictex \fi\fi

\def\gt{{\mathsurround=0pt\it $\cal G\mskip-2mu$eometry \&\ 
$\cal T\!\!$opology}}        

\def\gtp{{\mathsurround=0pt\it $\cal G\mskip-2mu$eometry \&\ 
$\cal T\!\!$opology $\cal P\!$ublications}}  

\def\lognumber#1{\def\thelognumber{#1}}
\def\volumenumber#1{\def\thevolumenumber{#1}}
\def\papernumber#1{\def\thepapernumber{#1}}
\def\volumeyear#1{\def\thevolumeyear{#1}}

\def\pagenumbers#1#2{\def\startpage{#1}\def\finishpage{#2}}
\def\published#1{\def\publishdate{#1}}
\def\proposed#1{\def\theproposer{#1}}
\def\seconded#1{\def\theseconders{#1}}
\def\received#1{\def\receiveddate{#1}}
\def\revised#1{\def\reviseddate{#1}}
\def\accepted#1{\def\accepteddate{#1}}
\def\asciititle#1{\def\theasciititle{#1}}

\def\asciiaddress#1{\def\theasciiaddress{#1}}
\def\asciiemail#1{\def\theasciiemail{#1}}
\long\def\asciiabstract#1{\long\def\theasciiabstract{#1}}

\let\\\par\let\thelognumber\relax
\let\thevolumenumber\relax\let\thepapernumber\relax
\let\thevolumeyear\relax\let\thesamplenumber\relax\let\startpage\relax
\let\finishpage\relax\let\publishdate\relax\let\receiveddate\relax
\let\reviseddate\relax\let\accepteddate\relax\let\theasciititle\relax
\let\theasciiauthors\relax\let\theasciiaddress\relax
\let\theasciiabstract\relax
\let\theasciiemail\relax\let\theshortauthors\relax\let\theshorttitle\relax

\long\def\maketitlep{   

\count0=\startpage

\gt\hfill      
\beginpicture
\setcoordinatesystem units <0.33truein, 0.33truein> point at 2.2 0.9
\setplotsymbol ({$\cal G$})
\plotsymbolspacing=9truept
\circulararc 315 degrees from 0 1 center at 0 0
\setplotsymbol ({$\cal T$})
\circulararc 315 degrees from 1 -1 center at 1 0
\endpicture
\break
{\small\ifx\thesamplenumber\relax 
Volume \else Sample
\fi\thevolumenumber\ (\thevolumeyear)
\startpage--\finishpage\nl
Published: \publishdate}
\vglue 0.5truein plus 0.4fil minus 0.1truein

{\parskip=0pt\leftskip 0pt plus 1fil\def\\{\par\smallskip}{\ifplaintex\large
\else\Large\fi\bf\thetitle}\par\medskip}   

\vglue 0pt plus 0.1fil 

{\parskip=0pt\leftskip 0pt plus 1fil\def\\{\par}{\sc\theauthors}
\par\medskip}

\vglue 0pt plus 0.1fil 

{\small\parskip=0pt\let\newline\\
{\leftskip 0pt plus 1fil\def\\{\par}{\sl\theaddress}\par}
\expandafter\ifx\theemail\relax    
\relax\else\vglue 5pt plus 0.02fil minus 2pt\def\\{\stdspace{\rm 
and}\stdspace} 
\cl{Email:\stdspace\tt\theemail}\fi
\ifx\theurl\relax                  
\relax\else\vglue 5pt plus 0.02fil minus 2pt\def\\{\stdspace{\rm 
and}\stdspace}
\cl{URL:\stdspace\tt\theurl}\fi\par}

\vglue 7pt plus 0.3fil minus 3pt

{\bf Abstract}
\vglue 5pt plus 0.1fil minus 2pt

\theabstract

\vglue 7pt plus 0.3fil minus 3pt

{\bf AMS Classification numbers}\quad Primary:\quad \theprimaryclass

Secondary:\quad \thesecondaryclass

\vglue 5pt plus 0.3fil minus 2pt

{\bf Keywords}\quad \thekeywords

\vglue 10pt plus 0.5fil minus 5pt

{\small  Proposed: \theproposer\hfill Received: \receiveddate\nl
Seconded: \theseconders\hfill 
\ifx\reviseddate\relax                         
Accepted: \accepteddate                        
\else
Revised: \reviseddate                          
\fi}
\eject
}       

\let\maketitlepage\maketitlep
\let\maketitle\maketitlepage

\font\phead=cmsl9 scaled 950
\font=cmsl9 scaled 1050
\font\pnum=cmbx10 scaled 913
\font\lnum=cmbx10 
\font\pfoot=cmsl9 scaled 950
\font=cmsl9 scaled 1050
\ifplaintex
\headline{\vbox to 0pt{\vskip -4.5mm\line{\small\phead\ifnum
\count0=\startpage ISSN 1364-0380 (on line)
1465-3060 (printed) \hfill {\pnum\folio}\else\ifodd\count0\def\\{ }%
\ifx\theshorttitle\relax\thetitle\else\theshorttitle\fi\hfill{\pnum\folio}
\else\def\\{ and }{\pnum\folio}\hfill\ifx\theshortauthors\relax\theauthors
\else\theshortauthors\fi\fi\fi}\vss}}
\footline{\vbox to 0pt{\vglue 0mm\line{\small\pfoot\ifnum\count0=\startpage
\copyright\ \gtp\hfill\else
\gt, Volume \thevolumenumber\ (\thevolumeyear)\hfill\fi}\vss
}}
\else
\makeatletter
\def\@oddhead{{\small\lhead\ifnum\count0=\startpage ISSN 1364-0380 (on line)
1465-3060 (printed) \hfill {\lnum\number\count0}\else\ifodd\count0
\def\\{ }\ifx\theshorttitle\relax \thetitle \else\theshorttitle\fi\hfill
{\lnum\number\count0}\else\def\\{ and }{\lnum\number\count0}
\hfill\ifx\theshortauthors\relax 
\theauthors\else\theshortauthors\fi\fi\fi}}\def\@evenhead{\@oddhead}
\def\@oddfoot{\small\lfoot\ifnum\count0=\startpage\copyright\ \gtp\hfill\else
\gt, Volume \thevolumenumber\ (\thevolumeyear)\hfill\fi}
\def\@evenfoot{\@oddfoot}
\makeatother
\fi

\newwrite\gtoutfile
\long\gdef\makeheadfile{  
{\def\\{, }\def\s{ }
\immediate\openout\gtoutfile head.xxx
\immediate\write\gtoutfile{To: math@arxiv.org}
\immediate\write\gtoutfile{Subject: put or rep NNNNN:pppp}
\immediate\write\gtoutfile{--text follows this line--}
\immediate\write\gtoutfile{Proxy-for: \ifx\theasciiauthors\relax
\theauthors\else\theasciiauthors\fi\s<\ifx\theasciiemail\relax\theemail\else\theasciiemail\fi>}
\immediate\write\gtoutfile{\noexpand\\}
\immediate\write\gtoutfile{Authors: \ifx\theasciiauthors\relax
\theauthors\else\theasciiauthors\fi}
\immediate\write\gtoutfile{Title: \ifx\theasciititle\relax
\thetitle\else\theasciititle\fi}
\immediate\write\gtoutfile{Subj-class: GT or SG or MG etc}
\immediate\write\gtoutfile{MSC-class: \theprimaryclass\ifx\thesecondaryclass\relax\else, \thesecondaryclass\fi}
\immediate\write\gtoutfile{Journal-ref: Geom. Topol. \thevolumenumber
(\thevolumeyear) \startpage-\finishpage}
\immediate\write\gtoutfile{Comments: Published by Geometry and Topology at}
\immediate\write\gtoutfile{\s\s http://www.maths.warwick.ac.uk/gt/GTVol\thevolumenumber/paper\thepapernumber.abs.html}
\immediate\write\gtoutfile{\noexpand\\}
\immediate\write\gtoutfile{}
\ifx\theasciiabstract\relax
\immediate\write\gtoutfile{\theabstract}\else
\immediate\write\gtoutfile{\theasciiabstract}\fi
\immediate\write\gtoutfile{}
\immediate\write\gtoutfile{\noexpand\\}
\immediate\write\gtoutfile{}
\immediate\closeout\gtoutfile}}  

\def\maketitlepage{\maketitlep\makeheadfile}
\let\maketitle\maketitlepage

\lognumber{161}
\volumenumber{5}\papernumber{25}\volumeyear{2001}
\pagenumbers{799}{830}
\received{12 January 2001}
\accepted{9 November 2001}
\proposed{Gang Tian}
\seconded{Yasha Eliashberg, Tomasz Mrowka}
\revised{9 October 2001}
\published{9 November 2001}

\usepackage{amssymb} 

\def\S{Section }

\let\oldsubsection\subsection

\def\subsection#1{\vspace{-10pt} \oldsubsection{#1}\vspace{-5pt}}
\def\nonumbersubsection#1{\vspace{-10pt} \oldsubsection*{#1}\vspace{-5pt}}

\newcommand{\R}{{\mathbb R}}

\newcommand{\Symp}{{\rm Symp}} 
\newcommand{\id}{{\rm id}}

\newcommand{\Hh}{{\mathcal H}} 
 
\newcommand{\Z}{{\mathbb Z}} 
\newcommand{\C}{{\mathbb C}} 
\newcommand{\oU}{{\overline {U}}}

\newcommand{\ou}{{\overline {u}}}
\newcommand{\ov}{\overline }
\newcommand{\oMm}{{\overline {\Mm}}}
\newcommand{\Ww}{{\mathcal {W}}}
\newcommand{\Xx}{{\mathcal {X}}}
\newcommand{\Yy}{{\mathcal {Y}}}

\newcommand{\p}{{\partial}} 
\newcommand{\al}{{\alpha}} 
 
\newcommand{\be}{{\beta}} 
 
\newcommand{\Om}{{\Omega}} 
\newcommand{\om}{{\omega}}
\newcommand{\io}{{\iota}}
\newcommand{\eps}{{\varepsilon}} 
\newcommand{\de}{{\delta}} 
 
\newcommand{\ga}{{\gamma}} 
\newcommand{\Ga}{{\Gamma}} 
 
\newcommand{\la}{{\lambda}} 
\newcommand{\si}{{\sigma}}

\newcommand{\Bb}{{\mathcal B}} 
\newcommand{\Cc}{{\mathcal C}} 
 
\newcommand{\Ll}{{\mathcal L}} 
\newcommand{\Jj}{{\mathcal J}} 
 
\newcommand{\Ff}{{\mathcal F}} 
 
\newcommand{\Nn}{{\mathcal N}} 
\newcommand{\Mm}{{\mathcal M}} 
\newcommand{\Pp}{{\mathcal P}} 
 
\newcommand{\Uu}{{\mathcal U}}

\newcommand{\La}{{\Lambda}} 
 
\newcommand{\grad}{{\rm grad\,}} 
\newcommand{\pr}{{\rm pr}} 

\newcommand{\Si}{{\Sigma}} 
\newcommand{\Ham}{{\rm Ham}} 
\newcommand{\area}{{\rm area\,}}

\newcommand{\MS}{{\medskip}}

\newcommand{\NI}{{\noindent}} 

\newcommand{\QED}{\endproof}

\newcommand{\CP}{{\rm CP}}

\newcommand{\cz}{{c_{HZ}}} 
\newcommand{\cs}{{c_{f}}} 
\newcommand{\Diff}{{\rm Diff}}

\newtheorem{theorem}{Theorem}[section] 
\newtheorem{cor}[theorem]{Corollary} 
 
\newtheorem{thm}[theorem]{Theorem} 
\newtheorem{prop}[theorem]{Proposition} 
 
\newtheorem{lemma}[theorem]{Lemma} 

\theoremstyle{definition}
\newtheorem{defn}[theorem]{Definition} 
\newtheorem{definition}[theorem]{Definition}

\newtheorem{rmk}[theorem]{Remark}
\newtheorem{guess}[theorem]{Conjecture} 
  
\begin{document} 

\title{Hofer--Zehnder capacity and length minimizing\\Hamiltonian paths}
\asciititle{Hofer-Zehnder capacity and length minimizing Hamiltonian paths}
\author{Dusa McDuff\\Jennifer Slimowitz}

\address{Department of Mathematics, State University of New York\\
Stony Brook, NY 11794-3651, USA}

\email{dusa@math.sunysb.edu}

\secondaddress{Department of Mathematics - MS 136, Rice University\\ 
Houston, TX 77005, USA} 
 
\secondemail{jslimow@rice.edu}

\asciiaddress{Department of Mathematics, State University of New York\\
Stony Brook, NY 11794-3651, USA\\and\\Department of Mathematics 
- MS 136, Rice University\\ 
Houston, TX 77005, USA}

\asciiemail{dusa@math.sunysb.edu\\jslimow@rice.edu}

\begin{abstract}
We use the criteria of  
Lalonde and McDuff to show that a path that is generated by a
generic autonomous
Hamiltonian is length minimizing with respect to the Hofer norm among all 
homotopic paths provided that
it  induces  no non-constant closed
trajectories in $M$. This generalizes a result of Hofer for 
symplectomorphisms of Euclidean space.  The proof for general $M$ uses 
 Liu--Tian's construction of $S^{1}$--invariant virtual moduli 
cycles.  As a corollary, we find that any semifree action of $S^{1}$
on $M$ gives rise to a nontrivial element in the fundamental group of the 
symplectomorphism group of $M$.  We also establish a version of the 
area-capacity inequality for quasicylinders.
\end{abstract}
 
\asciiabstract{We use the criteria of Lalonde and McDuff to show that a
path that is generated by a generic autonomous Hamiltonian is length
minimizing with respect to the Hofer norm among all homotopic paths
provided that it induces no non-constant closed trajectories in
M.  This generalizes a result of Hofer for symplectomorphisms of
Euclidean space.  The proof for general M uses Liu-Tian's
construction of S^1-invariant virtual moduli cycles.  As a
corollary, we find that any semifree action of S^1 on M gives
rise to a nontrivial element in the fundamental group of the
symplectomorphism group of M.  We also establish a version of the
area-capacity inequality for quasicylinders.}

\primaryclass{57R17}
\secondaryclass{57R57, 53D05}

\keywords{Symplectic geometry, Hamiltonian diffeomorphisms, Hofer norm, 
 Hofer--Zehnder capacity} 

\maketitle  
 
\section{Introduction} 
 
In this paper we provide a sufficient condition for a path
$\phi_t$, $0 \leq t \leq 1$  in the Hamiltonian
group $\Ham(M)$ to be length minimizing  with
respect to the  Hofer norm among homotopic paths with fixed endpoints. 
This extends the work done 
by Hofer~\cite{H}, Bialy--Polterovich~\cite{BP}, Ustilovsky~\cite{U},  and
Lalonde--McDuff~\cite{LM2} on characterizing geodesics in $\Ham(M)$.
We will work throughout on a closed symplectic manifold $(M, \om)$, 
though our results extend without difficulty to the group 
$\Ham^{c}(M, \om)$ of compactly supported Hamiltonian 
symplectomorphisms when $M$ is noncompact and without boundary.

There are several more or less equivalent definitions of the Hofer 
norm.  We will use Hofer's original definition.  Namely, we define the 
length $L(H_{t})$ of a  time dependent   Hamiltonian function
$H_t \co M \rightarrow {\bf R}$  for $0 \leq t \leq \tau$  to be  
$$
L(H_{t}) = \int_0 ^\tau \left(\max_{x \in M} H_t(x) - \min_{x
\in M} H_t(x)\right)\, dt.
$$ 
The length of the corresponding path
 $\phi_t^{H}, 0 \leq t \leq \tau,$ in $\Ham(M)$
is also taken to be $L(H_{t})$, and the Hofer norm $\|\phi \|$ of  
$\phi \in \Ham(M)$ is 
the infimum of the lengths of all of the paths from 
the identity to $\phi$.\footnote
{We fix signs by choosing $\phi_t^{H}$ to be tangent to the vector field 
$X$ defined by $\om(X, \cdot) = dH.$}
 This norm 
does not change if we restrict attention to paths parametrized by $t \in [0,1]$,
since this amounts to replacing $H_{t}, t\in [0,\tau],$ by
$\tau H_{\tau t}, t\in [0,1]$.  Hence, unless explicit mention is made to 
the contrary, all paths will be assumed to be so parametrized.

Although the Hofer 
norm is simply defined, it is difficult to calculate in general. 
One can separate this question into two: the first is to calculate the 
minimum of the lengths of paths between $\id$ and  $\phi$ in some fixed 
homotopy class, and the other is to minimize over the set of all homotopy 
classes.  We call paths that realise the first minimum {\it length 
minimizing in their homotopy class} (or simply {\it length minimizing}), 
and those realising the second
{\it absolutely length minimizing}.  It is  hard to find
absolutely length minimizing paths except in the very rare cases 
when $\pi_1(\Ham(M))$ is known.  However the first problem is often 
more manageable.  Also,  in cases where there  is a natural path 
from the identity to $\phi$ --- for example if there
is a path induced by   a circle action such as a rotation --- one can 
look for conditions under which this natural path is 
length minimizing. 
 
A simple  example of an absolutely 
length minimizing path is rotation of $S^{2}$ through $\pi$ radians: 
see \cite{LM2}, II  Lemma 1.7.   The proof
can be generalized to rotations of $\CP^2$ and of the one-point blow up of
$\CP^2$: see Slimowitz~\cite{Sl}. Because the argument uses explicit
embeddings of balls, it is too clumsy to work for general manifolds.

\subsection{Statement of main results}

In this note we concentrate on paths
 $\phi_t^H, t\in [0,1]$, that are generated by autonomous (ie time 
 independent) Hamiltonian 
 functions $H\co M\to \R$.  Our aim is to understand the
set 
$$
\La_H = \{\la: \mbox{the flow } \phi_t^{\la H}, t\in [0,1], \mbox{ of } \la H  
\mbox{ is length minimizing\hspace{2cm}}$$\vglue-35pt
$$\mbox{\hspace{8cm}in its homotopy class}\}.
$$
It is easy to see that this set is always a closed interval.  It  
has nonempty interior by Proposition~1.14 in~\cite{LM2}.\footnote
{The papers~\cite{LM2} were written at a time when 
it was not yet understood how to
define Gromov--Witten invariants for general 
symplectic manifolds $M$. Therefore,
many of the results in part II have unnecessary restrictions. 
    In particular, 
in Theorems 1.3 (i) and 1.4 and in Propositions 1.14 and 1.19 (i)
one can remove the hypothesis that
$M$ has dimension $\le 4$ or is semi-monotone.  The point is that 
these results rely on Proposition 4.1, and so use the fact that
 quasicylinders $Q = (M\times D^{2}, 
\Om)$ have the nonsqueezing property.   This is now known to hold for 
all $M$.}  
There are  Hamiltonians $H$ on manifolds with infinite 
fundamental group such that $\La_H = [0, \infty)$, ie, 
the flow of  $\la H$ is absolutely length minimizing for 
all $\la > 0$.  Here the lower bound for the length
is provided by the energy--capacity inequality on the universal cover:
see~\cite{LM2} Lemma~5.7.  When $M$ is closed and simply connected, 
in all known  examples (other than circle actions that have
 $\phi_1^{\la H} = id$ for some $\la >0$) 
the distance between the identity and
the symplectomorphism $\phi_{1}^{\la H}$ tends to infinity as 
$\la\to \infty$.\footnote{Added Dec 01:  In fact there are many other 
paths $\phi_1^{\la H}, \la \ge 0,$ that remain a bounded distance from 
$id$.  For example, if $F$ has support in a ball $B$ and $\psi (B)\cap B 
= \emptyset$, define $H = F - F\circ\psi$.  Then $\phi_1^{\la H} = 
\phi_1^{\la F}\circ\psi\circ (\phi_1^{\la F})^{-1}\circ\psi^{-1}$ 
remains at a distance $2\|\psi\|$ from $id$.}
However,  this path does not 
remain length minimizing for all $\la$.  Thus 
in the simply connected case 
one expects $\La_H$ to be a compact interval $[0, \la_{\max}(H)]$ 
for all $H$. 

The next result applies to all symplectic manifolds, and follows by an 
easy application of the curve shortening technique of~\cite{LM2} I 
Proposition 2.2.    

\begin{lemma}  Suppose that $H$ is a Hamiltonian  that assumes its 
maximum values on the set $X_{\max}$.  Then, if there is a Hamiltonian
symplectomorphism $\phi$ of $M$ such that $\phi(X_{\max})\cap X_{\max}$ 
is empty, $\la_{\max} < \infty$.
\end{lemma}

In this case one can estimate $\la_{\max}$ by comparing the 
displacement energy of a neighborhood $\Nn$ of $X_{\max}$ with the 
growth of $H$ on $\Nn$.
For a discussion of related questions see 
Polterovich~\cite{P}.

\MS

If $H$ is generic and hence a Morse function, it
follows from the above lemma that $\la_{\max} (H) <\infty$.
 However, one can get a sharper 
estimate for $\la_{\max}$ by looking at 
 the linearized flow near a critical point $p$ of $\la H$.  In suitable 
 coordinates, this has the form
$e^{- t\la JQ}$ where $Q$ is the Hessian of $H$ at $p$ and $J$ is the 
standard
almost complex structure. 
  We will say that $p$ is {\it overtwisted} for $\la H$
 if $A= -JQ$ has an imaginary  eigenvalue $i\mu$ with $ \mu > 2\pi$.  
This is equivalent to saying that the linearized 
flow of $\la H$ at $p$ has a nonconstant periodic orbit of 
period $< 1$: see \S\ref{ss:trans}.
  Ustilovsky's analysis in~\cite{U} of the second variation 
equation for geodesics shows that the path 
 $\phi_t^{\la H}, t\in [0, 1],$ ceases to be length minimizing as soon as all 
 the  global maxima  of $\la H$ are  overtwisted.
A similar result applies to minima, and also to certain 
degenerate $H$: see~\cite{LM2}. 

If $p$ 
 is an overtwisted local extremum of $H$, a celebrated result of 
Weinstein~\cite{We} implies that the nonlinear 
flow of $\la H$ near $p$ also has  nonconstant periodic orbits of 
period $< 1$.
Hence it is natural to make the following conjecture.

\begin{guess}\rm The path $\phi_t^H, 
t\in [0,1],$ is length minimizing in its homotopy class 
whenever its flow has no nonconstant 
contractible periodic orbits of period $< 1$.
\end{guess}

Hofer showed in~\cite{H} that this is true
for compactly supported Hamiltonians on 
$\R^{2n}$ by using a variational argument that 
does not extend to arbitrary manifolds: see also Section 5.7 in~\cite{HZ}. 
It was also established in the cases when
 $M$ has dimension two or is weakly exact
 in~\cite{LM2}~Theorem~5.4.   In this paper we extend the arguments
 in~\cite{LM2} to arbitrary manifolds.  Unfortunately this does not
 quite allow us to prove the full conjecture.
The problem is that there are 
functions $H$ with no nonconstant periodic orbits but yet with overtwisted 
critical points, and, for technical reasons, our argument cannot cope with
such points.   However, it is well 
known that for generic $H$ this problem does not occur; generic
overtwisted critical points always give rise to 
$1$--parameter families of contractible periodic orbits of period $< 1$.   For 
the sake of completeness, we 
give a simple topological proof of this in Lemma~\ref{le:nond} below 
and also describe Moser's example of an overtwisted Hamiltonian whose 
only periodic orbit is constant.

In view of this, it is useful to make the following definition.

\begin{defn}\label{def:non} \rm A periodic orbit is called {\it fast}
if its period is $< 1$.  Given an (autonomous) Hamiltonian $H$ we denote 
by $\Pp(H)$ the set of its fast contractible periodic orbits,
and by $\Pp_{crit}(H)$ the set of fast periodic orbits of the linearized 
flows at its critical points.  We will say that $H$ is {\it slow}
if the only elements in $\Pp(H)$ and $\Pp_{crit}(H)$ are 
constant paths.
\end{defn}

Here is the main  result of this paper.   
 
\begin{thm}\label{thm:main}
Given a  closed symplectic manifold $(M, \omega)$, let $\phi^H_t$,
$0 \leq t \leq 1,$ be the path in $\Ham (M)$  generated by 
the autonomous Hamiltonian  
$H\co M \rightarrow {\bf R}$.  If $H$ is slow,
 then  this path
is length minimizing among all homotopic paths between the identity and
$\phi^H_1$.   \end{thm}

Note that the path remains length minimizing in its homotopy class 
even if $H$ has periodic orbits of period exactly equal to $1$.  
To see this, first
apply the theorem to $(1 - \eps) H$ for  $\eps > 0$ 
and then use the fact that  the set $\La_{H}$ defined above is closed.

This theorem applies in particular to semi-free Hamiltonian circle actions
$\phi_t^H$, $t\in S^1 = \R/\Z$.  Recall that these are actions in which
the  stabilizer subgroups of each point are either trivial or the full 
group.  Thus in this case all nonfixed points lie on periodic orbits of period 
exactly $1$.  Moreover, because the flow $\phi_t^H$ on $M$ is congugate to its 
linearization near the critical points, it is easy to see that none 
of these points are overtwisted.

\begin{cor} Every semi-free symplectic $S^1$ action on a closed symplectic 
manifold $(M, \om)$ represents a nontrivial element $\ga$ 
in $\pi_1(\Symp(M,\om))$.  Moreover, if the action is Hamiltonian,
the corresponding loop has minimal length among all freely homotopic
loops in $\Ham(M,\om)$.
\end{cor}
\proof If the action is not Hamiltonian then the result is obvious 
(and the semi-free condition is not needed) since in this case 
the image of the loop 
under the flux homomorphism
$$
\pi_1\Symp(M, \om) \to H^1(M, \R)
$$
is nonzero.  For Hamiltonian loops, Theorem~\ref{thm:main} implies that
they are length minimizing paths from $id$ to $id$ in their homotopy class.
Because  the constant path to $id$ is always shorter 
than the given loop 
the latter cannot be null homotopic.  The last statement is an easy
consequence of the conjugacy invariance of the norm.
\QED
 
Somewhat surprisingly, there seems to be no elementary proof of
the first statement in
this corollary.  It would be interesting to know if it remains true in the 
smooth category.  In particular, do arbitrary smooth semi-free $S^1$ actions 
on $M$ represent nontrivial elements 
in $\pi_1(\Diff(M))$ or even in $\pi_1(H(M))$,
where $H(M)$ is the group of self-homotopy equivalences of $M$?  This is 
true for nonHamiltonian symplectic loops, since the flux homomorphism 
extends to $\pi_1(H(M))$.\footnote{Added in Dec 01: Claude LeBrun 
pointed out that the diagonal $S^1$ action on $\C^2$ given by multiplication 
by $e^{i\theta}$ induces a semifree action on $S^4$ that represents the 
trivial loop in $\pi_1(S0(5))\subset \pi_1(\Diff(S^4))$.  For further 
work on this subject see~\cite{MT}.}

Observe also that the semi-free condition is needed.  Consider, for example, 
the $S^{1}$ action on $\CP^{2}$ given by:
$$
[z_{0}: z_{1}: z_{2}]\mapsto 
[e^{i\theta}z_{0}: e^{-i\theta}z_{1}: z_{2}].
$$
This is null-homotopic, while 
 points such as $[1:1:0]$ have
 $\Z/2\Z$ stabilizer.  Clearly,  a general 
Hamiltonian $S^1$ action remains length minimizing for time $1/k$
where $k$ is the order of the largest isotropy group.

 As a byproduct of the proof we also calculate a very slightly 
 modified version of the Hofer--Zehnder capacity for cylinders $Z(a)$,
where
$$
Z(a) = (M \times 
D(a), \,\om \times \si_a)
$$
and $(D(a),\si_a)$ is a  $2$--disc 
 with total area $a$.   To explain this, we recall
the definition\footnote
{Hofer originally considered Hamiltonian systems in $\R^{2n}$ and 
hence had no need to restrict to contractible periodic orbits in 
condition (d) below.  In the definition of
$\cz$ given in~\cite{HZ}, this condition is not imposed.
We have inserted it here to make 
$\cz$ as relevant to our problem as possible.
This definition appears in Lu~\cite{Lu}, who pointed out
that the monotonicity axiom has to be suitably modified.
It is called the $\pi_{1}$--sensitive
Hofer--Zehnder capacity in  Schwarz~\cite{S}.}
 of the Hofer--Zehnder capacity $\cz$:
$$
{\cz}(N, \omega) = \sup \{ \max(H) \; \vline \; H \in {\cal H}_{ad}(N,
\omega) \}
$$  
where the set $ {\cal H}_{ad}(N, \omega)$ of admissible Hamiltonians
consists of all of the autonomous Hamiltonians on $N$  such that
\begin{itemize}
\item[\rm(a)]
For some compact set $K \subset N - \partial N$, 
$H|_{N - K} = \max(H)$ is constant;

\item[\rm(b)] There is a nonempty open set $U$ depending on $H$ such that  
$H|_ U = 0$;

\item[\rm(c)] $0 \leq H(x) \leq \max(H)$ for all $x \in N$;

\item[\rm(d)]  All  fast contractible periodic 
solutions of the Hamiltonian system $\dot x = 
X_H(x)$ on $N$ are constant. 
\end{itemize}

As explained above, our arguments are sensitive to the presence of
overtwisted critical points.  Hence we define the modified capacity
$c_{HZ}'$ as follows:
$$
c_{HZ}'(N, \omega) = \sup \{ \max(H) \; \vline \; H \in {\cal H}_{ad}'(N,
\omega) \}
$$  
where the set ${\cal H}_{ad}'(N, \omega)$ of admissible Hamiltonians
consists of all  autonomous Hamiltonians on $N$  
that satisfy conditions (a), (b), (c) above as well as the 
following version of (d): 
\begin{itemize}
\item[\rm(d$'$)]  $H$ is slow. 
\end{itemize}

These capacities are closely related. Clearly $c_{HZ} \le c'_{HZ}$.
Our discussion above
implies that
the set ${\cal H}_{ad}'(N, \omega)$ has second category in
${\cal H}_{ad}(N, \omega)$: see Corollary~\ref{cor:ot}.  Furthermore 
the two capacities may agree: it is not hard to see that they both equal $a$
on the $2$--disc $(D(a), \si_{a})$.\footnote
{In fact, there are no known examples where they differ.}
Since the capacity of the product $Z(a)$ is at least as large as 
that of $(D(a), \si_{a})$, the difficult part of the next 
proposition is to find an upper bound for $c_{HZ}'(Z(a))$.

\begin{prop}\label{prop:cylcap}
Let $(M,\om)$ be any closed symplectic manifold.
Then 
$$
c_{HZ}'(M\times D(a),\, \om\times \si_{a}) = a.
$$
\end{prop}

There are several ways in which one could try to generalize the main 
theorem.  Siburg showed in~\cite{Si} that the conjecture holds for
flows generated by time dependent Hamiltonians on ${\R^{2n}}$ provided
that these also have isolated and fixed extremal points.  (The fixed
extrema are needed to ensure that the path is a geodesic:
see~\cite{BP}.)  Although it seems very likely that
Theorem~\ref{thm:main} should hold on general $M$ in the time
dependent case, the method used here is not well adapted to tackle
this problem.  In fact, while our paper was being finished, Entov 
developed in~\cite{En} a rather different approach as
part of a larger program that has some very interesting applications. It 
may well be that his method would be better in the
time dependent case: see Remark~\ref{rmk:tdep}.

It is also natural to wonder what happens when $H$ does have 
nonconstant fast periodic orbits and/or overtwisted critical points.
 For example we might take an $H$ that satisfies the conditions of the 
 theorem and consider
the flow of $\la H$ for $\la > 1$.  It would seem plausible that if some 
critical point of index lying strictly between $0, 2n$ becomes overtwisted 
$\la H$ would remain length minimizing, at least for a while.
One problem here is that a critical point that is just on the point of 
beoming overtwisted (ie, has eigenvalue $2\pi i$) 
is degenerate as far as Floer theory is concerned.  The main step in
our proof is to demonstrate that a particular moduli space of Floer 
trajectories is nonempty, which we do by a deformation argument.
Thus we need to know that the relevant spaces of Floer trajectories 
are regular when $\la$ varies from $0$ to $1$, and it is here that 
the overtwisted critical points would cause a problem: see 
Lemma~\ref{le:crit}.
If degenerations occur, one must either carry through a 
detailed analysis of the degeneration or argue that this moduli space 
is nonempty for cohomological reasons.  Since both approaches would take us 
rather far from the main theme of this paper, we will not pursue them 
further here.

 \subsection{Techniques of proof}\label{ss:tech}

The proofs of the above results employ the criteria for length 
minimizing  paths developed in~\cite{LM2}.  For the convenience of the reader, 
this  is explained in \S2 below.  The idea is to compare the length of
the path with the capacity of an associated region in $M\times \R^{2}$ 
that is roughly speaking a cylinder.
In order to make the method 
work, it would suffice to  know that the 
Hofer--Zehnder capacity $\cz$ satisfies 
the {\it area--capacity inequality}
$$
\cz(Z(a)) \le a,
$$
on all cylinders.
This is equivalent to saying that every 
Hamiltonian $H\co Z(a) \to [0,c]$, that is identically zero on some open 
subset and equals its maximum value $c$ on a neighborhood of
the boundary $\p Z(a)$, has 
fast periodic orbits as soon as $c > a$. In~\cite{HV}, Hofer and Viterbo  
prove this statement for weakly exact $(M, \om)$, ie, when 
$\om|_{\pi_2(M)} = 0$.    Their argument was  extended to all manifolds by 
Liu--Tian in~\cite{LT}.  As these authors  point out,
 the ``usual'' theory of $J$--holomorphic  curves is not much help
 even in the semi-positive case
 because one must use moduli spaces on which there is an action of 
 $S^{1}$. Their paper establishes the needed technical 
basis --- $S^{1}$--equivariant Gromov--Witten invariants and virtual 
moduli cycles ---  to prove Proposition~\ref{prop:cylcap}
stated above.   However, they do not 
consider arbitrary Hamiltonians but a special class that is relevant 
to the Weinstein conjecture, and their paper is organised in such a 
way that one cannot simply quote the needed results.
This question is discussed further in \S\ref{ss:virt}.
 
In fact the above area--capacity  inequality  is
more than is needed for the problem at hand,
and  it is convenient to consider another modification of $\cz$ 
defined by maximizing over a restricted class
 of Hamiltonians that are compatible with the fibered structure of the 
cylinder.  This makes the geometry of the problem more 
 transparent and hence allows us to work with semi-positive
$M$ without using virtual 
 moduli cycles at all.
 
Here is a version of our main technical 
result.  It is somewhat simplified since we in fact need 
an analogous result to hold for quasicylinders, rather than just for 
cylinders: see \S2.
 It will be convenient to 
 think of the base disc $D(a)$ of $Z(a) = M\times D(a)$ as being a disc 
 on the Riemann sphere $S^{2} = \C\cup \{\infty\}$ with center at 
 $\infty$, and hence to call the central fiber
$M_\infty = M\times \{\infty\}$.

\begin{prop}\label{prop:HZ}  Let $F\co Z(a)\to [0,c],$ be a Hamiltonian function 
such that 
\begin{itemize}
\item[\rm(i)]  its only critical points occur in the sets
$M_\infty$ and $M\times \oU_0$, where
$\oU_0$ is a connected  neighborhood of the 
boundary  $\p D(a)$;

\item[\rm(ii)] near the central fiber $M_\infty$, $F = H_M + \be(r)$ where $H_M$ is
 a Morse function on $M$, and $\be$ is a function of the
  radial coordinate $r$ that is $< \pi r^{2}$ near $r=0$; 

\item[\rm(iii)]  $F\co Z(a)\to [0,c]$ is surjective, and is constant and 
equal to its maximum 
value on $M\times \oU_0$.
\end{itemize}
Then, if $c > a$, $F$ is not slow, ie, it  has either a nonconstant
fast periodic orbit or an overtwisted critical point.
\end{prop}

This paper is organized in the following way.   
The second section  describes the criteria for length minimizing 
paths developed by Lalonde and McDuff in \cite{LM2} and explains the 
role of  Hofer--Zehnder capacities.  
The third gives the proofs of the area--capacity inequality and  of 
  Proposition~\ref{prop:HZ}.   We discuss in detail some 
  technicalities about the intersections of bubbles and Floer 
  trajectories, that are omitted from standard references such 
  as~\cite{FHS}.
\MS

\NI {\bf Acknowledgements}\qua This paper is a development of part of
the second author's thesis.  The authors thank Helmut Hofer, Francois
Lalonde, GuangCun Lu, Leonid Polterovich, and Dietmar Salamon for very
helpful comments, and also Karen Uhlenbeck who pointed out a
significant gap in a much earlier version of the argument.  The first
author thanks Harvard University for providing a congenial atmosphere
in which to work on this paper.  The first author is partially
supported by NSF grants DMS 9704825 and 0072512.  The second author
was supported in 1998--9 by a grant awarded by the North Atlantic
Treaty Organization.

\section{Criteria for length minimizing paths} 
 
We  briefly describe the  Lalonde--McDuff 
 criterion for finding paths that are length minimizing in their homotopy
class.  In \cite{LM2}, they first  derive a geometric way of detecting that
$L(H_t) \leq L(K_t)$ for two Hamiltonians $H_t$ and  $K_t$ on $M$.  Then,
they determine sufficient   conditions involving symplectic capacities for this
geometric requirement  to be satisfied.   

For technical reasons it is convenient to restrict to Hamiltonians 
$H_{t}$ that are identically $0$ for $t$ near $0,1$.  This 
restriction does not cause any problems:  it is easy to 
see that every time independent Hamiltonian $H$ may be replaced by
one of the form $\be(t) H$ that satisfies the above condition 
and has the same length and time $1$--map as before.

\subsection{Estimating Hofer length via quasicylinders} \label{ss:outl}

 To begin, we must make a few definitions and set some notation. 
Suppose we have $H_{t}$, a time dependent Hamiltonian 
function on the closed symplectic manifold $(M^{2n}, \omega)$.  We may 
assume\footnote
{There is a slight technical problem here when the function
$
t\mapsto \min(t)= \min_{x \in M} H_t(x)
$
is not smooth.  In this case,
we replace $H_{t}$ by $H_{t} + m(t)$ where
$m(t)$ is a smooth function that is everywhere $\le \min(t)$
and is such that $\min(t) - m(t)$ has arbitrarily small integral.
This slightly changes the areas of the regions $R_{H}^{\pm}$.  
However, this can be absorbed into the $\nu$ fudge factor: 
we only need to measure lengths exactly for time 
independent  $H$.}
that  
for each $t$, 
$$
\min_{x \in M} H_t(x) = 0.$$  
We denote the graph $\Ga_H$ of $H_{t}$ by
$$
\Gamma_H = \{ ( x, H_t(x), t) \} \subset M \times \R \times [0,1].
$$  
Now, given some small $\nu > 0$  choose a function
$\ell(t)\co [0,1] \rightarrow [-2\nu, 0]$ such that
$ \int_0^1 - \ell(t) dt = \nu$.
A {\it thickening of the region under} $\Ga_H$ is 
$$
R_H^-(\nu) = \{ (x,s,t) \; \vline \; \ell(t) \leq s \leq H_t(x) \} 
\subset M \times [\ell(t), \infty) \times [0,1].
$$
Since $H_{t}\equiv 0$
for $t$ near $0,1$ we may arrange 
that $R_{H}^{-}$ is a manifold with corners along $s=0, t= 0,1$ by choosing 
the function $\ell(t)$ so that its graph is tangent to the lines $t=0, 
t=1$.

 Similarly, we can define 
$R_H^+(\nu)$ to be a slight thickening of the region above 
$\Ga_{H}$: 
$$
R_H^+({\nu}) = \{ (x,s,t) \; \vline \; H_t(x) \leq s \leq \mu_H(t) \} 
\subset M \times \R \times [0,1]
$$ 
where $\mu_H(t)$ is chosen so that   
$$
\mu_H(t) \geq \max_(t) = \max_{x \in M} H_t(x) ,\qquad \int_0^1 (\mu_H(t) -
\max_{t}) dt = \nu.
$$
  We define 
$$
R_H(2\nu)\; = \; R_H^-(\nu) \cup R_H^+(\nu)
\;  \subset 
M \times \R \times [0,1].
$$ 
 We equip $R_H^-(\nu)$, $R_H^+(\nu)$, and  
$R_H(2\nu)$ with the product symplectic 
form $\Omega = \omega \times \si$ where $\si = ds \wedge dt$. 
In particular, for any Hamiltonian $H_{t}$, $(R_H(\nu), \Omega)$ 
is symplectomorphic to  the product
$(M \times D(a), \Om)$ where 
$D(a)$ denotes the $2$--disc $D^2$ with area $a = L(H) + 2\nu$. 

Now, suppose $H_t$ and $K_t$ are two Hamiltonians on $M$ such that  
$\phi_1^H = \phi_1^K$ and the path $\phi_t^H$ for $ 0 \leq t \leq 1$  is 
homotopic (with fixed endpoints) to the path $\phi_t^K$ in $\Ham(M)$.   
There is a map $g\co\Gamma_K$ to $\Gamma_H$ defined by  
$$
g(x,s,t) = (\phi_t^H \circ (\phi_t^K)^{-1}(x), s - K(x)+H(\phi_t^H \circ 
(\phi_t^K)^{-1}(x)), t).
$$  This map $g$ extends to a symplectomorphism of  
$R_K^+(\nu)$, and we define  
$$
(R_{H,K}(2\nu), \Omega) =
R_H^-(\nu) \cup_g R_K^+ (\nu).
$$  
We assume that the functions $\ell$ and $\mu_H$ are chosen so that 
$R_{H,K}(2\nu)$ is a smooth manifold with boundary.  The 
contractibility of the  loop
$\phi_t^H \circ (\phi_t^K)^{-1}$ in $\Ham(M,\om)$ implies
that $(R_{H,K}(2\nu), \Omega)$  is
diffeomorphic to a product $(M\times D, \Om)$  by a diffeomorphism 
that is the identity near the boundary and is symplectic on each fiber.
However $\Om$ may not be a product, and so we make the following definition.

\begin{definition} \label{defn:qc} \rm
Let $(M, \omega)$ be a closed symplectic manifold and $D$ a set diffeomorphic 
to a disc in $({\bf R^2}, \si) $ where $\si = ds \wedge dt$.  Then, the manifold $Q =  
(M \times D, \Omega)$ endowed with the symplectic form $\Omega$ is called 
a {\it quasicylinder} if  
\begin{itemize} 
\item[\rm(i)] $\Omega$ restricts to $\omega$ on each fibre $M \times  
\{pt\}$; 
\item[\rm(ii)] $\Omega$ is the product $\omega \times \si$ near the
boundary  $M \times \partial D$.
\end{itemize} 
If $\Om  = \omega \times \si$ everywhere, not just 
near the boundary, $Q$ is called a {\it split} quasicylinder.
The {\it area} of any  quasicylinder $(M \times D, \Omega)$ 
is defined to 
be the number  $A$ such that 
$$ 
\mbox{ vol } (M \times D, \Omega) = A \cdot \mbox{ vol } (M, \omega).$$ 
Thus if $ (M \times D(a), \Omega)$ is split, its area is simply $a$.   
\end{definition}

Since  $(R_{H,K}(2\nu), \Omega)$  has trivial monodromy round its 
boundary, it is not hard to see that it is a quasicylinder:
see~\cite{LM2}~II\S2.1.  However, it may not be  split.

The key to the analysis is the following lemma
taken from~\cite{LM2}~II, Lemma 2.1, whose proof we include 
for the convenience of the reader.  It shows that
if the areas of both quasicylinders 
$(R_{H,K}(2\nu), \Omega)$ and $R_{K,H}(2\nu), \Omega)$ are greater than or 
equal to $L(H_t)$ for all $\nu$, then $L(H_t) \leq L(K_t)$.

\begin{lemma} 
\label{le:key} 
Suppose that $L(K_t) < L(H_t) = A.$  Then, for sufficiently small $\nu >0$, 
at least one of the quasicylinders $(R_{H,K}(2\nu), \Omega)$ and 
$(R_{K,H}(2\nu), \Omega)$ has area $< A$.   
\end{lemma} 
 
\proof 
Choose $\nu > 0$ so that  
$$L(K_t)+4 \nu < L(H_t).
$$ 
  Evidently, 
$$ \begin{array}{rcl} 
{\rm vol} (R_{H,K}(2\nu)) + {\rm vol} (R_{K,H} (2\nu)) & 
= & {\rm vol} (R_H(2\nu)) + {\rm vol} (R_K(2\nu)) \\ 
& = & ({\rm vol} M) \cdot (L(H_t) + L(K_t)+ 4\nu) \\ 
& < & 2({\rm vol} M) \cdot L(H_t)  
\end{array} $$ 
where $R_H(2\nu) = R_H^- (\nu) \cup R_H^+ (\nu)$.  \QED 

To 
proceed, one needs some way of finding lower bounds for the area of
a quasicylinder.  The arguments in~\cite{LM2}  use symplectic 
capacities, which are  functions from the set of symplectic manifolds to  
${\R} \cup \{ \infty \}$ satisfying certain properties; in particular, 
they are invariant under  symplectomorphisms.  

Suppose we have chosen a  
particular capacity $c$ and symplectic manifold $(M, \omega)$.  We say the   
{\it area--capacity inequality} holds for $c$ on $M$ if  
$$
c(M \times D, \Omega) \leq \mbox{ area of } (M \times D, \Omega)
$$ 
holds for all quasicylinders $(M \times D, \Omega)$.   It is 
useful to make the following definition.

\begin{definition}\rm The {\it capacity} $c(H_{t})$ of a Hamiltonian function 
$H_t$ is defined as 
$$
c(H_{t}) = \min \{ \inf_{\nu > 0} c(R_H^-(\nu)), 
\inf_{\nu > 0} c(R_H^+(\nu)) \}.
$$ 
\end{definition}   
 
Now, take a manifold $M$ and a capacity $c$ such that the area--capacity 
inequality holds for $c$ on $M$, and suppose that we have a Hamiltonian  
$H_t \co M \rightarrow { \R}$ for which 
$$
c(H_{t}) \geq L(H_t).$$ 
Then, for any Hamiltonian $K_t$ generating a flow $\phi_t^K$ which is homotopic 
with fixed end points to $\phi_t^H$ (and thus has $\phi_1^K = \phi_1^H$), 
we can embed $R_H^-(\nu)$ into $R_{H,K}(2\nu)$ and
$R_H^+(\nu)$  into $R_{K,H}(2\nu)$.  Thus, we know 
$$ 
L(H_t) \;\leq\; c(H_{t}) \;\leq \;c(R_H^-(\nu)) \;\leq\; c(R_{H,K}(2\nu))
$$ 
$$ 
L(H_t)\; \leq\;c(H_{t}) \;\leq c(R_H^+(\nu)) \;\leq \;c(R_{K,H}(2\nu)),
$$ 
with the last inequality in both lines holding by the 
monotonicity property of capacities.  Since the area--capacity inequality holds,
we know that the  areas of both quasicylinders $R_{H,K}(2\nu)$ and
$R_{K,H}(2\nu)$ must be greater than or equal to their capacities  and hence
greater than or equal to $L(H_t)$.  Therefore, by Lemma~\ref{le:key},  
$L(K_t) \geq L(H_t).$  This proves the following result (Proposition~2.2 from
\cite{LM2}, Part II.) 
 
\begin{prop}\label{prop:crit}  Let $M$ be any  
symplectic manifold and $H_{t \in [0,1]}$ a Hamiltonian 
generating an isotopy $\phi_t^H$ from the identity to $\phi = \phi_1^H$.  Suppose 
there exists a capacity $c$ such that the following two conditions hold: 
\begin{itemize} 
\item[\rm(i)] $c(H_{t}) \geq L(H_t)$ and 
\item[\rm(ii)] for all Hamiltonian isotopies $\phi_t^K$ homotopic 
rel endpoints to $\phi_t^H$ , $t \in [0,1]$, the  
area--capacity inequality holds (with respect to the given capacity $c$) for
the  quasicylinders $R_{H,K}(2\nu)$ and $R_{K,H}(2\nu)$.  
\end{itemize} 
Then,  the path $\{\phi_t^H\}_{t\in [0,1]}$ minimizes length 
among all homotopic
Hamiltonian paths from $id$ to $\phi$.
\end{prop} 
 
Hence, to show that $H_t$ generates a length minimizing 
 path $\{\phi_t^H\}_{t \in [0,1]}$, 
we need only produce a capacity $c$ that satisfies the above conditions (i) 
and (ii).  Various results were obtained in~\cite{LM2} by using
the Gromov capacity $c_{G}$ and the 
Hofer--Zehnder capacity $\cz$.
It seems to be best to use $\cz$, since condition (i) holds for it
almost by definition whenever $H$
has no nontrivial fast periodic orbits, while
(i) is very restrictive for $c_{G}$.  On the other hand, the existence of 
 Gromov--Witten invariants on general symplectic manifolds allows one 
 to show  easily that condition  (ii) holds for $c_{G}$, while 
the proof of (ii) for $\cz$ is more subtle.  Liu--Tian consider a 
very closely related question in~\cite{LT}, and using their methods 
one can prove that (ii) holds for the very slightly modified version
$c_{HZ}'$ of $\cz$ on any manifold: 
see~\S\ref{ss:virt}.  
 
 In view of the complexity of the constructions in~\cite{LT}, we  present in 
 the next section a different modification of
the Hofer--Zehnder capacity for which one can prove 
 condition (ii) without too much difficulty in the semi-positive case.
  This capacity $\cs$
 is defined  for fibered spaces such as quasicylinders, 
 satisfies (i) whenever $H$ is slow and also satisfies (ii) for any closed $M$. 
It depends on some extra structure that we
need to choose and so is not defined for all symplectic manifolds.  
Note that the only
properties of the capacity $c$ that we used above are that it is defined for sets such as
$R_H^\pm(\nu)$ and that it has the monotonity property $$
c(R_H^-(\nu)) \le c(R_{H,K}(2\nu)),\qquad
c(R_H^+(\nu)) \le c(R_{K,H}(2\nu)).
$$

\subsection{The Hofer--Zehnder capacity for fibered spaces} \label{ss:modHZ}
 
We first explain what is meant by a fibered symplectic manifold.

\begin{defn}\label{defn:fib}  \rm  We will say that the symplectic manifold $(Q, \Om)$
is {\it fibered} with fiber $(M, \om)$ if there is a submersion $\pi\co Q\to D^2$
such that $\Om$ restricts to a nondegenerate form on each fiber $M_b =
\pi^{-1}(b)$, where $(M_{b}, \om_{b})$ is symplectomorphic to $(M, 
\om)$ for one and hence all $b$.  In this case, because $D^{2}$ is 
contractible one can use Moser's theorem to choose an identification 
$s_{Q}$ of $Q$ with $M\times D^{2}$ so that $\om_{b} = \om$ for all $b$.
$s_Q$ is said to  {\it normalize}  $Q$ if in addition there is a
 small closed disc  $\oU_\infty$  in $D^2$ with center $\infty$ so 
 that $\Om$ restricts to $\om \times\si$ on $M\times U_{\infty}$,
 where $\si$ is the area form $ds\wedge dt$ as before.  A symplectic embedding
$\psi\co Q\to Q'$ is said to be {\it  normalized}
if it takes the central fiber $M_\infty$ in $Q$ to that in $Q'$ and if
$$
\psi = (s_{Q'})^{-1} \circ s_Q
$$
on some neighborhood of $M_\infty$ that need not be the whole of
$\pi^{-1}\oU_\infty$. 
\end{defn}

Using the symplectic neighborhood theorem 
it is easy to see that every fibered space
can be normalized near any fiber.  
Further, every  quasicylinder $(Q, \Om)$ is fibered, though in general
the identification $Q \to M\times D^2$
that occurs in the definition of a quasicylinder 
is a normalization only near fibers that
are sufficiently close to the boundary. It is also not hard to see that 
the spaces $(R_H^\pm(\nu), \Om)$ can be fibered with
fibers  $\pi^{-1}(b)$ of the form $\{(x, s_b(x), t_b): x\in M\}$: the restriction of $\Om$
to such sets equals $\om$ since $t_{b}$ is fixed.  We will assume that the
fibers lying in the part of $R_H^-(\nu)$ with $s < 0$ are flat, ie, 
also have fixed
$s$--coordinate $s_b(x) = s_{b}$.  This normalizes $R_H^-(\nu)$  near some fiber
$M_0$ with $s < 0$.  Similarly,  the fibration of $R_H^+(\nu)$ is chosen 
to have flat
fibers $s = const$  near its upper boundary $ s = \mu_H(t)$.  This 
means that spaces such as $R_{H,K}(\nu)$ have two possible 
normalizations, one at a fiber where $s<0$ and the other near its 
upper boundary.  However, it is not hard to see that there is a 
fiberwise symplectomorphism taking one to the other so that they are 
equivalent.

\begin{defn}\label{def:adf}\rm\hspace{-4.3pt}
Given a normalized fibered space $Q$, we define the set 
$\Hh_{f,ad}(Q)$ of admissible
Hamiltonians  to be the set of all functions $F\co Q\to [0,\infty) $ such that:
\begin{itemize}
\item[\rm(i)] in some neighborhood $M\times \oU_\infty$ of 
the central fiber $M_{\infty}$, 
$F = H_M + \be(r)$ where $H_M$ is
 a Morse function on $M$, and $\be$ is a function of the radial coordinate
 $r$ of the disc that is $< \pi r^{2}$;
 
\item[\rm(ii)]  $F\ge 0$ everywhere and is constant and equal to its 
maximum on a product neighborhood 
$M\times \oU_{0}$ of the boundary;

\item[\rm(iii)]  the only critical points of $F$ occur on $M_{\infty}$
and in $M\times \oU_{0}$;

\item[\rm(iv)]  $F$ is slow. 
\end{itemize}
\end{defn}

\begin{defn}\label{def:cHZmod}  \rm  We define the  Hofer--Zehnder capacity
of a normalized fibered space $Q$ by
$$
\cs(Q) = 
\sup \{ \max(F) \; \vline \; F \in {\cal H}_{f,ad}(Q)\}
$$  
\end{defn}

Clearly,  this capacity $\cs$ has the appropriate monotonicity property, ie,
$\cs (Q) \le \cs(Q')$ whenever there is a normalized symplectic embedding $Q \to Q'$.
In particular,
$$
\cs(R_H^-(\nu)) \le \cs(R_{H,K}(2\nu)),\qquad 
\cs(R_H^+(\nu)) \le \cs(R_{K,H}(2\nu)).
$$
 The following proposition, which is proved in \S3, shows that $\cs$ also
satisfies
condition (ii) in Proposition~\ref{prop:crit}.

\begin{prop}\label{prop:HZ1}
For any normalized quasicylinder $(Q, \Om)$ of area $A$,
$$
\cs(Q) \le A.
$$
\end{prop}

We next check condition (i). 

\begin{lemma} \label{le:autH} If $H\co  M\to \R$ is slow, then $\cs(H) \ge L(H)$.
\end{lemma}
\proof  This is essentially~\cite{LM2} II, Proposition 3.1.  We will prove 
that $\cs(R_H^-(\nu))$ $\ge L(H)$.  The case of $R_H^+(\nu)$ is similar: 
indeed $R_H^+(\nu)$ is symplectomorphic to
 $R_{m-H}^-(\nu)$, where $m = \max H$.

By assumption, $H$  has minimum value $0$.  Let $m$ be its maximum, and 
consider the set  
$$
S_{H,\nu} = \{(x,\rho,\tau)\in M\times D(m+\nu/2)\; | \: 0 \le \rho \le H(x) + \nu/2\},
$$
where $(\rho, \tau)$ are the action-angle coordinates on the disc given in terms of
polar coordinates $(r, \theta)$ by
$$
\rho = \pi r^2,\quad \tau = \frac \theta {2\pi}.
$$
This space $S_{H, \nu}$ is essentially the same as $R_H^-(\nu)$.
Indeed, it is not hard to check that there is a symplectic  embedding  
 $S_{H, \nu} \to R_H^-(\nu)$ of the form $(x, \rho, \tau)\mapsto (x, 
 \phi(\rho, \tau))$
for some area preserving map $\phi\co  \R^2\to \R^2.$  Moreover, $S_{H, \nu}$ is fibered
with central fiber at $(\rho,\tau) = (0,0)$,
and we may choose this embedding so that it respects suitable normalizations of
both spaces.   Hence it suffices to show that for all $\eps > 0$
$$
\cs(S_{H, \nu}) \ge L(H) - \eps.
$$

To see this, first consider the function  $F = m - H(x) +\rho$.  This 
is constant and equal to $m + \nu/2$ on $\p S_{H, \nu}$, and its flow is 
given by
$$
\phi_F^t\co  (x,\rho,\tau)\mapsto (\phi_H^t(x), \rho, \tau + t).
$$
Since $H$ is slow and the critical points
of $H$ give rise to  periodic orbits for $F$ with period precisely 
$1$, $F$ is also slow.  Now smooth out $F$ to $F_{\eps}:
S_{H, \nu} \to \R$, where
$$
 F_{\eps}(x,\rho, \tau) = 
\left\{\!\!\!\begin{array}{lll}
(1-\eps)\left(m -H(x) + \al_\nu(\rho)\right), & \mbox{if} & \rho < \nu/4,\\
(1-\eps)F(x,\rho,\tau) &  \mbox{if} & \nu/4\le  \rho \le H(x) + \nu/4,\\
(1-\eps)\left(m -\al_\nu(H(x) +\nu/4 - \rho)\right),\kern -4.4pt &
 \mbox{if}& H(x) + \nu/4 \le \rho \le\\ 
&&\hspace{1.8cm}H(x) +  \nu/2.
\end{array}\right.
$$
Here $\eps > 0$, and $\al_\nu(\la)$ is a
increasing smooth surjection $\la\co  [0,\nu] \to [0,\nu]$ 
that is $\le \la^{2}$ near $0$  and equals $\la$
when $\la \ge \nu/6$.   
Since the flow of $(1 - \eps)F$ goes slower 
than that of $F$ when  $\eps > 0$,  $(1 - \eps)F$ is slow.
Now the bump function $\al_\nu(\rho)$ must have
derivative slightly $> 1$ somewhere.  Hence when we turn it on the flow in the
$\tau$--direction goes slightly faster.  However, for each given $\eps$ 
we can clearly
choose $\al_\nu$ so that the product $(1-\eps)\al_\nu(\rho)$ is slow.
A similar remark applies to the smoothing at $\p S_{H,\nu}$.   
Hence $F_\eps$ is slow and has maximum value $m - \eps = L(H) - \eps$.    

If $H$ were a Morse function, 
 $F_\eps$ would be admissible, ie, belong to $\Hh_{f,ad}(S_{H, \nu})$, 
 and the proof
would be complete.
Hence the last step is to alter $F_\eps$ near the central fiber 
by replacing $H$ with a
function  that is independent of $\rho$ for $\rho $ near $0$ and 
restricts to a Morse function $H_{M}$ on $M_{\infty}$.
 This is easy to do without introducing any nonconstant fast 
 periodic orbits
 since we just need to change $H$ in directions along which 
its second derivative is
small.  See, for example, Lemma 12.27 in~\cite{MSintro} that 
shows that $H$ is slow whenever 
its second derivative is
sufficiently small.\QED

\NI
{\bf Proof of Theorem~\ref{thm:main}}

This follows by the preceding lemma and by the remarks at the end of 
\S\ref{ss:outl}.

\begin{rmk}\label{rmk:tdep}\rm
Suppose that $H_{t}$ is a time dependent Hamiltonian.  The space 
$R_{H}^{-}$ is again essentially the same as 
 $S_{H,\nu}$ where this is defined to be
 the set of points $(x,\rho,\tau)$ with $0 \le \rho\le
H_{\tau}(x)$, and we can define the (time independent)
Hamiltonian $F$ near its boundary $\p S$ to be (a 
smoothing of) $m - H_{\tau}(x) + \rho$ as before.
The problem is that this function is not well defined on the central 
fiber $M_{\infty}$ since $\tau$ is not a coordinate there, and there 
seems to be no satisfactory way of understanding when one can make such an
extension.  In particular, it seems 
one would need the restriction of $F$ to $M_{\infty}$ to 
have the same norm as $H_{t}$ and yet be slow.
Entov in~\cite{En} 
connects the Hamiltonian $H$ to the geometry of a fibered space via 
the choice of suitable connection rather than by the construction of 
the Hamiltonian $F$.  The condition on the connection 
is local while our condition on $F$ (that it should be slow) is global. 
Hence his approach seems better 
adapted to this problem.  
\end{rmk}

\section{The area--capacity inequality}\label{s:ca}

We begin by sketching the proof of this inequality 
for semi-positive $M$ using  
the set up in  Hofer--Viterbo~\cite{HV}.  \S\ref{ss:trans} contains 
more technical details, and \S\ref{ss:virt}  discusses the case of 
general $M$.

\subsection{Outline of the proof}\label{ss:sket}

For simplicity, we will assume for now
 that $M$ is semi-positive, ie, that  one of 
the following conditions holds:\MS

\NI
(a)  the restriction to $\pi_{2}(M)$ of the first Chern class 
$c_{1}(M)$ of $M$ is positively proportional to $[\om]$ -- the 
monotone case; or\MS

\NI
(b)   the minimal Chern number $N$ of $M$ is $> n-2$, where $2n = \dim 
M$. \MS

In this case the Gromov--Witten invariants on $M$ can be defined naively, ie, 
bubbles can be avoided, simply by choosing a generic $J$ on $M$: 
see~\cite{MS}.  It  is not necessary to use 
the virtual moduli cycle.  Notice that usually one asks that $N > 
n-3$ in (b).  Strengthening this requirement allows us to say that 
no element of a
generic $2$--parameter family of almost complex structures on $M$
admits a holomorphic curve of negative Chern number.

We will assume in what follows that
 $(Q,\Om)$ is a quasicylinder and that $F$ is an admissible 
 Hamiltonian in the sense of
Definition~\ref{def:adf}.  In particular, this means that for all 
$\la \le 1$ the only $1$--periodic orbits of the flow of $\la F$ on 
$M_{\infty}$ are constant and occur at the critical points $p_{k}$ of $F$.
Thus every Floer trajectory for $\la F$ on $M_{\infty}$ converges to 
these critical points.  Our aim is to show:

\begin{prop} If $F$ is an admissible Hamiltonian on 
     the quasicylinder\break $(Q,\Om)$ and if $M$ is semi-positive then
    $\|F\| \le \area Q$.
    \end{prop}

Because $(Q,\Om)$ is a product near its boundary 
 $\p Q$ we  can identify this to a single fiber $M_{0}$ and so replace $Q$ 
 by  $(V = M\times S^{2}, \Om)$ where  $\Om$ restricts to $\om$
 on each fiber.

\begin{defn}\rm  
 An $\Om$--tame almost complex 
 structure $J$ on $V$ will be said to be
{\it normalized} if each fiber is $J$--holomorphic
 and  if  in addition it is a 
product near both $M_0$ and $M_{\infty}$.  \end{defn}

Thus each such $J$ defines 
a $2$--parameter family of $\om$--tame almost complex structures 
on $M$, and by our assumptions on $M$ we can assume that there are no  
$J$--holomorphic spheres that have Chern number $< 0$ and lie in a 
fiber of $V$.  Since the existence of such curves is what necessitates 
the introduction of virtual moduli cycles, we will be able to count 
curves in $V$ (and hence define appropriate 
Gromov--Witten invariants) provided that we are in a situation 
where the only bubbles that 
appear lie in its fibers.

 The idea of the proof is to assume that $\|F\| > \area Q$ and to find a 
contradiction. 
Let $A  = [pt\times S^{2}]\in H_{2}(V)$.  
It is shown in~\cite{LM2} that there is a family of noncohomologous
symplectic forms  $\Om_{s}$ on $V$ starting with $\Om_{0} = \Om$ such 
that $\Om_{1}$ is a product.  Hence the fibered space  $(V, \Om)$ is
deformation equivalent to a product, which implies that 
 $Gr(A) = 1$, where the Gromov invariant $Gr(A)$ 
 counts the number of 
 $J$--holomorphic $A$--spheres in $V$ going 
through some fixed point $p$ in $V$ for sufficiently generic $J$.  
We will choose $p$ to be some minimum  $p_{\infty}\in M_{\infty}$ 
of $F$, and will fix the parametrizations $u$ of the spheres
by requiring that 
 $$
 u(0)\in M_{0},\;\; u(1)\in 
 M_{1},\;\; u(\infty) = p_{\infty}\in M_{\infty},
 $$
 where $M_{1}$ is some fiber distinct from $M_{0}, M_{\infty}$.
The arguments given in \S3.2 below show that one can calculate 
$Gr(A)$ using  generic normalized $J$.  Hence, for such $J$ the 
number of these curves will sum up to $1$ when  counted with the 
appropriate signs.  (In fact, in this semi-positive case, one can use 
mod $2$ invariants and so ignore the sign.)

We now ``turn on'' the perturbation 
corresponding to the Hamiltonian flow of
$\la F$ for increasing $\la \ge 0$.\footnote
{One must be very careful with signs here since there are many 
different conventions in use.  We have chosen to use the upward 
gradient flow of $F$ (even though it is more usual to use the downward 
flow) because this fits in with our set-up.  Since $F$ 
takes its maximum on $M_{0}$ we need to consider trajectories going from 
this maximum to a minimum: see Lemma~\ref{le:la} below.}
The resulting trajectories $u$
have domain $\C$ and in terms of the coordinates 
 $(s,t)$ of $(-\infty,\infty)\times S^{1}$ satisfy the following equation
 for some $\la$:
\begin{eqnarray}\label{eq}
\p_{s} u  + J(u) \p_{t}(u) &=&  \la\, (\grad F) \circ u,\\
\lim_{s\to-\infty} u(s,t) \in M_{0},&& 
\lim_{s\to \infty} u(s,t) = p_{\infty},
\end{eqnarray}
where $\grad F$ is the gradient of $F$ with respect to the metric 
defined by $\Om$ and $J$.
Because $dF = 0$ near $M_{0}$ the map $u$ is $J$--holomorphic for $s 
<<0$ and so, by the removable singularity theorem, does extend to a 
holomorphic map $\C\to V$.  Thus $u$ is a generalized Floer trajectory 
of the kind considered in ~\cite{HV,PSS}, and we will call it a 
$\la$--trajectory.
Because its limit at $\infty$ is a point, 
it also extends to a continuous map $S^{2}\to V$ that represents the 
class $A$.  It is shown in~\cite{HV} that the algebraic number of 
solutions to this equation is still $1$ for small $\la$.
 
Given $F$ and a normalized $J$, let $\Cc = \Cc_{A}$ be the moduli space 
consisting of all pairs $(u, \la)$ where $\la \in [0,1]$ and
$u\co  \R\times S^{1}\to V$ 
satisfies equations~(1), (2) as well as the following 
normalization condition:\MS

\NI
$(*) \;\;\;$  $u(0,0) \in M_{1}$ where $M_{1}$ is a fiber of $Q$ distinct 
from $M_0, M_{\infty}$.
\MS

Note that $\Om(A)$ 
is precisely the area of $Q$.  The crucial ingredient that ties the 
solutions of the above equation
 to the area--capacity inequality is the fact that the 
size $\|F\|$ of $F$ gives an upper bound for $\la$.

\begin{lemma}\label{le:la}  If $(u, \la) \in \Cc_{A}$ then $\la \|F\| < 
\Om(A)$ = area $Q$.
\end{lemma}
\proof A standard calculation shows that the action functional
$$
a(s) = \int_{(-\infty,s]\times S^{1}} u^{*}\Om + \int_{0}^{1} \la 
H(u(s,t)) dt
$$ 
is a strictly increasing function of $s$.
Since $F(p_{\infty}) = 0$ and $F|_{M_{0}} = \|F\|$
by construction, the action $a(s)$ satisfies
$$
\lim_{s\to-\infty} a(s) = \la \|F\|,\quad \lim_{s\to\infty} a(s)= \Om(A). 
$$
Hence $\la \|F\| < \Om(A)$ as claimed.  \QED

Note that if $p_{\infty}$ is a nonovertwisted critical point of $F$ of 
Morse index $k$, 
then the formal dimension of $\Cc$ is $1 + k$ (see for 
example~\cite{PSS}) and so equals $1$  with 
the current choice of $p_{\infty}$.
Because $A$ is not a multiple class, it follows from the standard 
theory  that
for any $M$  we can regularize the moduli space $\Cc$ by choosing a 
generic normalized $J$: see \S3.2.  Hence for such a choice $\Cc$ is 
a manifold of dimension $1$ lying over $[0, 1]$ via the projection 
$$
\pr\co  \Cc\to [0,1],\qquad (u,\la)\mapsto \la.
$$ 
 Because $\la$ is  restricted to the interval $[0,1]$, $\Cc$ could 
have boundary over $\la = 0,1$.
As mentioned above, $0$ is a 
regular value for $\pr$  for generic $J$,
and the algebraic number of points in $\pr^{-1}(0)$ 
is $1$.  On the other hand, we know from Lemma~\ref{le:la} above 
that, if $\|F\| \ge \area Q$, the set $\pr^{-1}(\la)$ is empty for $\la= 1$.
 The only way to reconcile these statements is for $\Cc$ to 
be noncompact.\MS

\nonumbersubsection{Noncompactness of $\Cc$}

Noncompactness in a moduli space of $J$--holomorphic Floer
trajectories is\break caused either by the bubbling off of $J$--holomorphic 
spheres or by the splitting of  Floer trajectories.  Now bubbling is a 
codimension $2$ phenomenon, and so, provided that we can make 
everything regular by choosing a suitably generic $J$, it will not occur 
along the $1$--dimensional space $\Cc$. It is easy to see that all 
bubbles have to lie in some fiber.  Hence, by our choice of 
normalization for $J$, we can avoid all bubbles.
(There are some extra details here that are discussed  in
\S\ref{ss:trans} below.)

Floer splitting  is harder to 
deal with since it occurs  in codimension $1$: a generic 
$1$--parameter family of Floer trajectories can degenerate into a pair 
of such trajectories.  For example, the trajectories in $\Cc$ could 
converge to the concatenation of a $\la$--trajectory $u\co \C \to V$ 
in class $A - B$ that converges 
to some critical point $p_{k}$ on $M_{\infty}$ of index $k$ together with a Floer 
$\la$--trajectory in $M_{\infty}$
from $p_{k}$ to $p_{\infty}$ in class $B\in H_{2}(M)$.  We will see 
in Lemma~\ref{le:fiber} below that these are the only 
degenerations that happen generically. 
Observe also that these degenerations do not occur in the situation treated 
by Hofer--Viterbo because of their topological assumptions on $M$.

To analyse this situation further,
denote by 
$$
\Cc_{A - B}(p_{k})
$$
the space of all pairs $(u, \la)$, 
where  $u\co \C\to V$ is a solution
to equations~(1), (2) with $p_{\infty}$ replaced by 
$p_{k}$, that
 is normalised by condition $(*)$ and represents the class $A-B$.
Similarly, denote by 
$$
\Ff = \Ff_{B}(p_{k})
$$
the space of all
pairs $(v, \la)$ where $v\co  \R\times S^{1} \to M_{\infty}$ is a Floer 
trajectory for $\la F$ from $p_{k}$ to $p_{\infty}$ in class $B$.
Note that the classes $B$ that occur here are constrained by the 
inequality $\om(B) < \om (A)$.  Moreover, since our assumption is that 
$\|F\| > \area Q$, we can slightly perturb $F$ within the class of 
admissible Hamiltonians to make $H_{M}$ slow and generic in the sense of
 Lemma~\ref{le:crit}.  
 That lemma then says that  we can choose  $J$
so that all the relevant moduli spaces of simple trajectories are 
regular, ie, have dimension equal to their formal dimension.
Thus $\Cc_{A - B}(p_{k})$  will have dimension 
$-2c_{1}(B) + k + 1$, where $k = {\rm index\,}  p_{k}$.
Further if  $B\ne 0$ is a simple (ie nonmultiple) class, then  ${\Ff}$ has
dimension $2c_{1}(B) - k + 1$. Because $F$ and 
$J_{M}$ are independent of the time coordinate $t$
and because the trajectories in $\Ff$ limit on fixed points rather 
than nonconstant periodic orbits, there is a
$2$--dimensional reparametrization group acting on the trajectories  
in $\Ff$.  Thus we need $2c_{1}(B) - k + 1\ge 2$ for $\Ff$ to be
nonempty, while we need $- 2c_{1}(B) + k + 1 \ge 0$ for 
$\Cc_{A - B}(p_{k})$ to be nonempty.  Therefore, if these spaces are 
both nonempty, $\Ff$ has dimension $2$ and $\Cc_{A - B}(p_{k})$ has 
dimension $0$.  Hence these spaces both consist of discrete sets 
of points, which, for generic $J$, will project to disjoint sets in the
$\la$--parameter space.  Thus  this kind of 
degeneration does not occur for generic $J$.

The crucial point in this argument is that the elements in $\Ff$
have an $S^{1}$ symmetry.   This presents a problem, since in general one cannot 
regularize  Floer moduli spaces containing multiply covered
trajectories unless one allows either the Hamiltonian $F$ or the almost complex 
structure $J$ to depend on $t$: see~\cite{FHS}.  The usual way to deal with 
this is to assume that $M$ is monotone: see Floer~\cite{F}.  However, we now 
show that in our special situation this assumption is unnecessary.

First observe that  we must also
avoid the case when the trajectory itself is independent of $t$, since 
then the $S^{1}$ action becomes vacuous.   But this could only happen 
if $B = 0$ and our
choice of $p_{\infty}$  implies both that 
$k\ge 0$ and that $B\ne 0$.  (Because the action $a(s)$ is 
strictly increasing and $F(p_{k}) \ge F(p_{\infty})$ we must have 
$\om(B) > 0$.)  The above argument shows that we need 
$2c_{1}(B) - k + 1\ge 2$ and hence $c_{1}(B) > 0$
for $\Ff$ to be nonempty when $B$ is simple and $J$ is generic. 
Moreover, if there is a multiply 
covered trajectory in class $\ell B, \ell > 1,$ from $p_{k}$ to 
$p_{\infty}$  then it covers an 
underlying simple trajectory in class $B$ between these points.
Therefore we must have
$c_{1}(B) > 0$ and $2c_{1}(B) - k + 1\ge 2$ in this case too.  But then the 
formal dimension
$- 2\ell c_{1}(B) + k + 1$ of $\Cc_{A - \ell B}(p_{k})$ is always 
negative.  But, because $A - \ell B$ is not a multiple class,
this moduli space consists of simple trajectories.  Therefore our assumptions 
imply that it is regular and hence empty for generic $J$.

It follows (modulo a few details discussed in \S\ref{ss:trans} below) that
there are no degenerations of the trajectories in $\Cc$ for $\la \in 
[0,1]$.  But we saw earlier that if $\|F\| \ge \area Q$
these trajectories must degenerate, ie, $\Cc$ cannot be 
compact.  Therefore
$\|F\| < \area Q$.

We have used the fact that 
none of the critical points of $F$ are overtwisted twice in the above 
argument.  First, it implies that the contribution of each critical point 
$p_k$ to the dimension of $\Cc_A$ is just its Morse index $k$ and so 
is $\ge 0$.  Second, we need the 
space of $\la$--trajectories to $p$ to be regular for {\it each} $\la \in [0,1]$
which is impossible if the 
linearized flow at $p$ has a periodic orbit of period $\la$.

\subsection{More details}\label{ss:trans}

We first discuss the behavior of the flow near overtwisted critical 
points, and then give more details of the transversality arguments 
needed to understand the compactification of $\Cc$.

\nonumbersubsection{Overtwisted critical points}

Since this question is local, we consider Hamiltonians $H\co  
\R^{2n}\to \R$ with a nondegenerate critical point at $0$.  We denote 
the Hessian by $Q$ so that the linearized flow at $0$ is $e^{At}$ 
where $A = -J_{0}Q$.  The eigenvalues of $A$ occur in real or imaginary
pairs $\pm \la, 
\pm i\la$, $\la\in \R$, or in quadruplets $\pm \mu, \pm \ov{\mu}, 
\mu \in \C - (\R \cup i\R)$.  Correspondingly, $\R^{2n}$ decomposes as a 
symplectically orthogonal sum of eigenspaces, one for each pair or 
quadruplet.  We will be concerned with the partial decomposition 
$$
\R^{2n} = E \oplus \sum_{j=1}^{k} E_{j}
$$ 
where the purely imaginary eigenvalues of $A$
are $\pm \la_1,\dots,\pm \la_k$ and
$E_{j}\otimes \C$ is the sum of the eigenspaces
for the pair
$\pm i\la_j$, and $E\times \C$ is the sum of the others. 
Observe that each $E_{j}$ contains a subspace of dimension at least 
$2$ that is filled out by periodic orbits of $e^{At}$ of period
$2\pi/\la_{j}$.  Indeed, for each eigenvector $v\in \C^{2n}$ in 
$E_{j}\otimes \C$ the intersection of $E_{j}$ with the subspace $\C\, v
\oplus \C \, \ov v$ consists entirely of such periodic orbits.  Hence, 
if $A$ has imaginary eigenvectors the linearized flow always has 
nonconstant periodic orbits. 

However this is not necessarily true for the nonlinear flow 
$\phi_{t}^{H}$.  
Moser considers the following example in~\cite{Mo}:\footnote
{He uses complex variables.  Observe that if $z_{k} = x_{k} + i y_{k}$
the Hamiltonian flow with our sign conventions 
can be written as $\dot z_{k} = -2i (\p H/\p 
{\ov z}_{k})$.}
$$
H(z_{1}, z_{2}) = \frac 12(|z_1|^2 - |z_2|^2) + (|z_1|^2 + |z_2|^2) 
\Re(z_{1}z_{2}).
$$
Clearly, the eigenvalues of $A_{H}$ are $\pm i$.  However, it is easy 
to check that the time derivative of the function $\Im(z_{1}z_{2})$ is 
strictly negative whenever $(z_{1}, z_{2})\ne (0,0).$  Hence there are 
no nonconstant periodic orbits.  

The problem here is that the two 
eigenvalues are equal.  More generally, similar phenomena can occur if any 
pair $i\la, i\la'$ of eigenvalues are resonant, ie, if  the ratio 
$\la'/\la$ is integral.
The next result is well known, and is proved in the real analytic case 
in Siegel--Moser~\cite{SM} \S16.  

 \begin{lemma}\label{le:nond}  Suppose in the above situation that
     $i\la$ is an imaginary eigenvalue of $A$ of 
     multiplicity $1$ that is  nonresonant in the sense that
     the ratio $\la'/\la$ is nonintegral for all other imaginary 
     eigenvalues $i\la'$ of $A$.  Then the flow $\phi_{t}^{H}$ of 
$H$ has a periodic orbit of period 
close to $2\pi/\la$ on every energy surface
close to zero.
 \end{lemma}
 \proof The linearized flow around $\{0\}$ is $e^{At}$ where
$A = -J_{0}Q $. As above $\R^{2n}$ decomposes as a
symplectically othogonal sum
$E_{0} \oplus E_{\la}$, where $E_{\la}$ is a $2$--dimensional space 
filled by periodic orbits of period $2\pi/\la$ and the restriction 
of $A$ to $E_{0}$ has no eigenvalues of the form $ik\la, k\in \Z$.
Consider the level set
$$
S_{1} = \{x\in \R^{2n}: H_{Q}(x) = 1\}
$$
of the quadratic part $H_{Q}$ of $H$.  By construction, 
it intersects $E_{\la}$ in a 
periodic orbit $\ga$ for $e^{tA}$  of period $T= 2\pi/\la$. 
The first return map 
$\phi_{\ga}$ of this orbit can be identified with the restriction
$e^{TA_{0}}$ of $e^{TA}$ to $E_{0}$.  Hence our assumptions on the 
eigenvalues of $A$ imply that its only fixed 
point is at the origin.  Thus its Gauss map 
$$
g\co  S^{2n-3}\to S^{2n-3},\quad v\mapsto \frac{\phi_{\ga}(v) - v}
{\|\phi_{\ga}(v) - v\|}
$$
is well defined.  Observe that $g$ has degree $1$.  In fact it is 
injective.
For, otherwise there would be vectors $v,w$ lying on different rays
in $E_{0}$ such that $\phi_{\ga} (v) - v = \phi_{\ga} (w) - w$.  
Since $\phi_{\ga}$ is linear, this would imply that it has $1$ as an 
eigenvalue, contrary to hypothesis.

Now consider the functions
$x\mapsto \eps^{-2}H(\eps x)$.  Since they
converge to $H_{Q}$ as $\eps$ decreases to $0$, for each fixed 
sufficiently small $\eps$
the orbits that start near $\ga$ remain near $\ga$ for $t\in [0,T]$. 
Hence the first return map given by following these orbits 
round $\ga$ is a perturbation $\phi_{\ga}^{\eps}$  of $\phi_{\ga}$. 
Hence its Gauss map is also defined and has degree $1$ for small 
$\eps$.  But this means that the Gauss map cannot extend over the 
interior of $S^{2n-3}$; in other words, $\phi_{\ga}^{\eps}$ must
have a fixed point.
This corresponds to 
 a closed periodic orbit of $\eps^{-2}H(\eps x)$ that is close to 
 $\ga$ and has period $T_{\eps}$ 
 close to $T$.  Since $\eps^{-2}H(\eps x)$ is conjugate to $H$,
 this implies that $H$ also has a periodic orbit of 
 period $T_{\eps}$.  
 \QED

\begin{cor}\label{cor:ot}
If a generic $H$ has an overtwisted critical point, ie, if its Hessian 
has imaginary eigenvalue $i\la$ with $\la > 2\pi$, then its flow has a 
nonconstant periodic orbit of period $< 1$.
\end{cor}
\proof  The hypotheses of the above lemma are satisfied for 
generic $H$.\QED

\begin{lemma}\label{le:crit}  
Suppose that the Hamiltonian $H$ on $M$ is slow.
  Then $H$ has arbitrarily small perturbations
$H'$ such that for generic $J$  the moduli spaces of simple Floer 
trajectories for $\la H', \la \in [0,1],$ in classes $B\in H_{2}(M)$
are all regular.
\end{lemma}
\proof If necessary, we first replace $H$ by $cH$ for some $c$ close to $1$ 
so that neither $H$ nor its linearized flows have nonconstant periodic
orbits of period $\le 1$.  Then slightly perturb $H$ so that it is also
a Morse function.  Finally, note that 
by~\cite{FHS} Remark 7.3
we may perturb $H$ to $H'$ so
that for all $\la\in [0,1]$
the critical points of  $\la H'$ satisfy 
the nondegeneracy conditions  of~\cite{FHS}~Lemma~7.2  with respect to 
a generic set of $J$ and for all $\la$.  
Thus simple (ie nonmultiply covered)
Floer trajectories all have regular injective points in the sense 
of~\cite{FHS}~\S7.   
The result now follows by
\cite{FHS}~Theorem 7.4. \QED

As always, it is not enough to know that trajectory spaces are 
regular.  One also needs to show that their closures have the right 
dimension.  This will follow from  Lemma~\ref{le:evaltr} below.
 
\nonumbersubsection{Structure of the stable maps in the closure of $\Cc$}
 
Next let us check that the degenerations of the elements in $\Cc$ 
really are compatible with the fibration.  By the standard compactness 
theorem, these degenerations consist of a finite number of 
Floer $\la$--trajectories $u_{i}\co \R\times S^{1}\to V$, $i = 
\ell,\ldots,k$ that are laid end to end 
together with some bubbles $v_{j}\co S^{2}\to V$.  Here, the $u_{i}$ 
are labelled in order, so that
$$
\begin{array}{llll}
\lim_{s\to \infty}  u_{i} &=& \lim_{s\to -\infty}  u_{i+1},&  \ell <i < k.
\end{array}
$$
Since the only critical points are either near $M_{0}$ or on $M_{\infty}$  there 
has to be at least one trajectory going between these manifolds. 
Pick one of them and call it $u_{1}$.  (We will see that in fact there is 
only one such trajectory.)
Because $F$ is slow,  the $u_{i}$ converge  
to critical points of $F$ at each end and so represent some homology 
classes in $V$.  In the proof of the next result it is convenient to 
allow ourselves to decrease the component $\be(r)$ of $F$ that is 
perpendicular to the fiber at $M_{\infty}$.  Since we assumed $\be < 
\pi r^{2}$ for small $r$,
we can reduce $\be$ to $\eps r^{2}$ on $r < \de/2$ for 
any $\eps$ without introducing any nonconstant fast periodic orbits. 

\begin{lemma}\label{le:fiber} Let $(u_{i}, v_{j})$ be a limit
of elements of $\Cc$ as described above.  If $\eps$ is sufficiently small,
 each bubble $v_{j}$  is contained in some fiber, and
the $u_{i}, i\ne 1,$ are Floer $\la$--trajectories in $M_{\infty}$.  
Moreover, $\ell = 1$ and the homology class represented by $u_{1}$ 
has the form $A - B$, for some $B\in H_{2}(M)$ with $0 \le \om (B) < \om (A).$  
\end{lemma}
\proof Suppose that $(u^{\al},\la^{\al})$ is a sequence of elements of 
$\Cc$ that converges weakly to a limit of the above type, 
where $u^{\al}\co  \C\to V$. 
Fix $\al$ and consider the composite map 
$$
\ou^{\al}= \pi\circ u^{\al}\co  \C \to V \to S^{2}.
$$
Since $J$ 
is a product near $M_{0}$ this map  is holomorphic over the inverse 
image of the neighborhood $\oU_{0}$ of $0\in S^{2}$.  Hence, because it has 
degree $1$, the projection from the 
image of $u^{\al}$ to the base is injective over $\oU_{0}$.

 Let $z_{j}$ be the set of points in $\C$ at which
$|du^{\al}(z)|\to \infty$. Then the restriction of $u^{\al}$ 
to compact pieces of $\C - \cup z_{j}$ converges to a map whose projection to 
the base is 
holomorphic and nonconstant over $\oU_{0}$.  Thus this limit is the 
trajectory $u_{1}$.   Since its intersection with the fiber class is 
$1$, it must represent some class of the form $A - B$, with $B\in 
H_{2}(M)$.  

Now consider the bubbles.  These are always $J$--holomorphic and so 
their projections to the base are holomorphic near $M_{0}$. 
Further, because the fibers are $J$--holomorphic they intesect each 
fiber positively.  Hence each bubble either is entirely contained
in a single fiber or represents
 a class $kA + B$ with $k> 0$.  But in the latter case they must intersect 
each fiber of
$M\times \oU_{0}$ which is impossible  because the projection from the 
image of $u^{\al}$ to the base is injective over $\oU_{0}$ and, 
as noted above, these 
points converge to the component $u_{1}$.  

Finally, consider the Floer trajectories.  Suppose there was a 
trajectory that came before $u_{1}$ and so 
had endpoint on $M_{0}$.
The previous argument applies to show that it is entirely contained in
$M_{0}$ and therefore satisfies the unperturbed Cauchy--Riemann 
equation and should be considered as a bubble.  In particular there is 
only one Floer trajectory that meets both $M_{0}$ and $M_{\infty}$
namely $u_{1}$. Hence the other Floer 
trajectories begin and end at points in $M_{\infty}$, and we claim that
for sufficiently small $\eps$ they are completely contained in 
$M_{\infty}$.

To see this, note that if $\eps$ were $0$, then $F$ would depend only 
on the fiber coordinates in the neighborhood  $r <\de/2$ of 
$M_{\infty}$.  Thus the Floer trajectories would project
to holomorphic trajectories in the base 
 and positivity of intersections with the fiber would imply
as before that 
the trajectories are entirely contained in $M_{\infty}$.  Therefore, 
because we are only interested in trajectories lying in a finite set 
of homology classes and with a finite set of possible endpoints, 
standard compactness arguments imply that for sufficiently small 
$\eps$ all trajectories must be contained in the neighborhood
$M_{\infty}\times \{r < \de/2\}$ of $M_{\infty}$.  Thus these 
trajectories would project to nullhomologous 
Floer trajectories in $S^{2}$ for the 
function $\eps r^{2}$ that begin and end at the point $r = 0$.  But 
these do not exist because the action functional could not increase 
strictly along such a trajectory.

It remains to prove the statement about the class $A-B$ represented 
by $u_{1}$.  Let $B_{i}, B_{j}$ be the classes represented by the other
$u_{i}$ and the bubbles $v_{j}$.  Clearly each $\om (B_{j}) > 0$.
Further each $\om(B_{i}) > 0$ because $a$ strictly increases along 
each trajectory and $p_{\infty}$ is a minimum of $F$: see 
Lemma~\ref{le:la}.  Similarly, $\om(A-B) > 0$ since $u_{1}(0)$ lies 
at a maximum of $F$.  Since
 $\om(B)$ is the sum of the $\om(B_{i}), \om(B_{j})$, the result 
 follows.  
\QED

\nonumbersubsection{Transversality of intersections of bubbles with trajectories}

First observe that by the previous lemma the only classes $B\in H_{2}(M)$ that 
occur as a component $u_{i}$ or $v_{j}$ of a limiting trajectory
in the closure of $\Cc$  have $\om (B) < \om(A) = \area Q$. Hence
only a finite number of classes can occur.  As already noted, standard theory 
tells us that we can regularize the moduli spaces of
vertical bubbles in $V$ and make all their 
intersections transverse by choosing generic normalized $J$ on $V$.
Thus all spaces of bubble trees (or cusp-curves) can be assumed 
to be of the right dimension.

Similarly, as we noted in Lemma~\ref{le:crit}, spaces of nonmultiply covered 
Floer trajectories in $M_{\infty}$ as well as the moduli spaces
$\Cc_{B,p_{k}}$ can be regularized 
by a time independent $J$ by~\cite{FHS}.  Thus there is a subset $\Jj_{reg}$ of 
second category in the space of all normalized almost complex 
structures on $Q$ such that all spaces of bubble trees and of 
simple trajectories are regular.

In order to make the ``usual'' theory of $J$--holomorphic curves work 
we must also ensure that these moduli 
spaces intersect transversally.  The basic arguments that establish 
this for spheres are given in~\cite{MS}  and the case of Floer 
trajectories is discussed in~\cite{FHS}.  
However, the standard proof that spaces of bubbles can be assumed to 
intersect transversally uses the fact that if
 two distinct simple bubbles ${\rm im\,} u$ and $ {\rm im\,} v$ intersect at
 some point $x = u(z) = v(w)$ then there is a small annulus $\al$ 
 around $z$ whose image by $u$ does not intersect ${\rm im\,} v$: see~\cite{MS} 
 Propositions 6.3.3 and 2.3.2.  This holds because otherwise the 
 two curves are infinitely tangent at $x$ and so must coincide.  This 
 argument breaks down for
  bubbles and Floer trajectories since they satisfy 
 different equations.   Since this detail seems to 
 have been ignored in  standard references such as~\cite{FHS}, we 
 deal with it now.

 For simplicity, we will suppose that there is just one bubble and so 
 will consider the intersection of the space of 
 unparametrized bubbles in class $B$ with the moduli space 
 $\Cc_{B'} =\Cc_{B',p_{\infty}}$.  It suffices to consider the 
 intersection of the corresponding parametrized curves.  Hence let
 $\Xx$ be the space of all maps 
 $$
 u\co (S^{2}, 0, 
 \infty)\to (Q, M_{0}, p_{\infty})
 $$
 in the class $A - B'$, let
 $\Yy$ be the space of all maps $v\co  S^{2}\to Q$ representing the 
 class $B$, and consider 
 the space $\Uu$ of all tuples
 $$
(u,v, \la, z, J)\in \Xx\times \Yy\times \R\times 
S^{2}\times \Jj
$$
satisfying the following conditions:
\begin{itemize}
\item[\rm(i)]  $u$ is a Floer $\la$--trajectory with respect to $J$;

\item[\rm(ii)] the bubble $v$ is $J$ holomorphic.
\end{itemize}
We want to show that when $J$ lies in 
a subset $ \Jj_{{reg}}$ of second category in $\Jj$ the space
$$
\{(u,v,z): (u,v,\la,z,J)\in \Uu, u(z) = v(0)\}
$$
is a manifold of the correct dimension.  This follows in the usual 
way from the next lemma.

\begin{lemma}\label{le:evaltr}  The evaluation map
$$
ev\co   \Uu \to Q\times Q: (u,v,\la, z, J)\mapsto (u(z),v(0))
$$
is transverse to the diagonal.
\end{lemma}
\proof  If $z = 0$ then $u$ is $J$--holomorphic near $z$ and the 
 argument of \cite{MS} Propositions 6.3.3 works.  The case $z = \infty$
is somewhat special since the moduli space of $u$--trajectories does not
have a tangent space at this point.  However, this does not matter since 
 $u(z)$ is fixed for all $J$ because it is the endpoint of the 
 Floer trajectory.  Instead we look at the space of $v$--bubbles
 and can appeal to Theorem~6.1.1 of 
 \cite{MS} that says that the map from the space of all pairs 
 $(v,J)$ in $\Uu$ to $Q$ given by evaluation 
 $$
 ev_{2}\co (v,J)\mapsto v(0)
 $$
 is surjective.

 When $z\ne 0, \infty,$ we can identify the domain of $u$ with 
 $\C$ and
 by reparametrization fix $z = 1$. The domain of the 
 linearization $D_{u}$ of the 
defining equation for the Floer trajectory equation  at $u$ is then
the space  $W^{1,p}(u^{*} TQ)$ which is defined to be the closure
with respect to the $(1,p)$--Sobolev norm of the space of compactly 
supported $C^{\infty}$--sections of $u^{*}TQ$ that are tangent to
the fiber at $z=0$: see~\cite{FHS} \S5.  Thus we may replace $\Uu$
by the space $\Uu'$ of tuples $(u,v,\la, J)$.  The
 tangent space of $\Uu'$ 
  at $(u,v,\la,J)$ consists of elements $(\xi_{1}, \xi_{2}, r, Y)$ 
  with $\xi_{1}\in W^{1,p}(u^{*}TQ)$, $\xi_{2}\in W^{1,p}(v^{*}TQ)$
  and  such that  
\begin{eqnarray*}
D_{u}(\xi_{1}) +\frac 12 Y(u)\circ du\circ i = rg_{F},& & (*)\\
D_{v}(\xi_{2}) +\frac 12 Y(v)\circ dv\circ i = 0 & & (**).
\end{eqnarray*}
(Here $g_{F}$ is the appropriate term coming from the variation in 
$\la F$.)
Moreover the derivative $d(ev)$ of the evaluation map is given by
$$
d(ev)(\xi_{1}, \xi_{2}, Y) = (\xi_{1}(1), \xi_{2}(0))\in 
T_{(x,x)}(Q\times Q).
$$
We know by Theorem~6.1.1 in~\cite{MS} that the map $(\xi_{2}, Y) \to 
\xi_{2}(0)\in T_{x}Q$ is surjective.  Hence given $a\in T_{x}Q$ 
there is $(\xi_{2}^{a}, Y^{a})$ that satisfy (**) with
$\xi_{2}^{a}(0) = a$.  Note that we cannot assume that the 
support of $Y^{a}$ is disjoint from the image of $u$ though we can 
make it in an arbitrarily small neighborhood of the intersection 
point $v(0)$.  Thus the 
element 
$
\nu = \frac 12 Y^{a}\circ du \circ i$ may well be nonzero.  Clearly, it 
will suffice  
to find $(\xi_{1},Y)$ so that
$$
\xi_{1}(1) = 0,\quad L(\xi_{1},Y) = -\nu,\quad Y=0 \mbox{ in the 
support of }Y^{a}
$$
where 
$$
L(\xi_{1},Y) = D_{u}(\xi_{1}) +\frac 12 Y(u)\circ du\circ i.
$$
The usual proof of transversality (as in~\cite{MS} Proposition 3.4.1
or ~\cite{FHS} Theorem 7.4) 
shows that the operator $L$ is surjective if $\xi_{1}$ ranges freely 
in $W^{1,p}(u^{*}TQ)$ and $Y$ is constrained to have support near any 
injective point of $u$.  In particular,  the condition that $\xi_1(0)$ be 
tangent  to the fiber can be fulfilled by adding a suitable vector 
tangent to the  group of M\"obius transformations of $S^2$ that fix 
$\infty$ and $1$.
Since the image of $v$ lies in a fiber distinct from $M_{0}$ and $u$ is 
injective near there we can easily arrange that the support of $Y$ 
is disjoint from that of $Y^{a}$.   Thus the only problem is the 
question of how to deal with the condition $\xi_{1}(1) = 0$.

To do this, we must consider more closely the proof that $L$ is 
surjective.  The argument goes as follows.  Since 
$$
D_{u}\co   W^{1,p}(u^{*}TQ) \to L^{p}(\Om^{0,1} u^{*}TQ)
$$
is Fredholm, 
the image of $L$ is closed and it suffices to show that it is dense.  
If not, there is $\eta$ in the dual space  $L^{q}((\Om^{0,1} 
u^{*}TQ)^{*})$
that vanishes on ${\rm im\,} L$.  In the standard case this implies 
that $\eta$ is a weak solution of the adjoint equation 
$D_{u}^{*}\eta = 0$ since it vanishes on all the elements 
$D_{u}\xi_{1}$.  Hence, by elliptic regularity, it is a strong 
solution of this equation.  It also must vanish in some open set because
it pairs to zero with all the elements $L(0,Y)$.  Hence $\eta = 0$ as 
required.

In our case $\xi_{1}$ is not an arbitrary element of
$W^{1,p}(u^{*}TQ)$ but rather is in the image of the map
$$
W^{1,p}(u^{*}TQ\otimes E) \stackrel{\phi}\to W^{1,p}(u^{*}TQ)
$$
where $E$ is a holomorphic bundle over $S^{2}$ with Chern class $-1$ and $\phi$
tensors the sections of $u^{*}TQ\otimes E$ by a holomorphic 
section $s$ of the dual bundle $E^{*}$ that vanishes at $1$.  Since 
$s$ is holomorphic
there is a commutative diagram
$$
\begin{array}{ccc}
W^{1,p}(u^{*}TQ\otimes E) & \stackrel{D_{u}^{E}}\to & 
L^{p}(\Om^{0,1} u^{*}TQ\otimes E)\\
{\otimes s} \downarrow  & & {\otimes s}\downarrow\\
W^{1,p}(u^{*}TQ) &\stackrel{D_{u}}\to & 
L^{p}(\Om^{0,1} u^{*}TQ).
\end{array}
$$
It follows that the image $\eta^{E} = \phi^{*}(\eta) = \eta\otimes s$ of 
$\eta$ in 
$L^{q}((\Om^{0,1} u^{*}TQ\otimes E)^{*})$ is a weak 
solution of the adjoint equation  
$(D_{u}^{E})^{*}\eta^{E} = 0$.  The standard argument applies to show 
that $\eta^{E}= \eta\otimes s$ is zero. Hence the $L^{q}$--section $\eta$ also
vanishes.\QED

\subsection{The case of general $M$}\label{ss:virt}
 
To construct the virtual moduli cycle as in~\cite{LT2} for curves in 
some manifold $(V, \om)$ one 
looks at the configuration space $\Bb$ of all pointed stable maps in some 
class $A$ that are 
nearly holomorphic.  Roughly speaking, $\Bb$ is an orbifold that
supports a orbibundle $\Ll$ whose fiber $L_{u}$ at the map $u\co \Si\to V$
is the Sobolev space of $L^{k,p}$--smooth sections of the bundle
$\Om^{0,1}(\Si, u^{*}(TV))$ of $(0,1)$--forms on the nodal Riemann surface 
$\Si$.  For each $J$, the delbar operator $\ov{\p}_{J}$ defines a section 
of $\Ll$ whose zero set is the set
 $\oMm_{J}$ of $J$--holomorphic stable maps.  If the derivative 
 $$
 D_{u}\co  L^{k+1,p}(\Si, u^{*}(TV)) \to L_{u}
 $$
 of this 
 map is surjective for all $(\Si, u)
 \in \oMm_{J}$, this zero set is an orbifold of the right dimension and its 
 fundamental cycle can be used 
 to define Gromov--Witten invariants.  Although  $\oMm_{J}$ 
is always compact with 
respect to the weak topology of $\Bb$, it might well
 be that for all $J'$ near $J$ this derivative is badly behaved, so that $\oMm_{J'}$ 
has components of too large dimension.  What one does to remedy 
the situation is define, 
 over some orbifold neighborhood $\Ww$ of
$\oMm_{J}$ in $\Bb$,  a finite-dimensional subspace $R$ of the set of 
sections of $\Ll$ 
such that the map
 $$
 D_{u} \oplus \io_{u}\co  L^{k+1,p}(\Si, u^{*}(TV))\oplus R \to L_{u}
 $$
is surjective for all $(\Si, u)$ in some smaller neighborhood $\Ww_{R}$
of $\oMm_{J}$, where $\io_{u}$ denotes evaluation 
at $u$.  This implies that for a generic 
small element $\nu\in R$ the set of solutions of  the perturbed 
Cauchy--Riemann equation 
$$
\p_{J}(u) + \io_{u}(\nu) = 0
$$
has the right dimension and supports a fundamental cycle.  
This is often called the virtual moduli cycle or regularized moduli 
space $\oMm\,\!^{\nu}$.

This is the briefest outline of Liu--Tian's method. Many more details can be 
found in~\cite{LT2,LT3,Mcv}.
The main point is the construction of $R$.  The idea is to find a 
suitable perturbation space $R_{i}$ over each subset $U_{i}$ of an open cover of 
$\oMm_{J}$ and then to patch these together.  

In our situation we start with an action of $S^{1}$ by reparametrization 
on the space of $J$-holomorphic Floer trajectories in $V=M$ 
between two points  $p$ and $q$ and want to construct the 
regularization $\oMm\,\!^{\nu}$
so that it also supports an $S^{1}$--action.  To do this one must 
first extend the original action to the neighborhood $\Ww$.  This 
extension 
will not  simply be an action of $S^{1}$: if a trajectory splits into 
two, or more generally $k$, pieces there will be an $S^{1}$ action on each part,
and one has to make everything equivariant with respect to this.  In 
particular, one must choose the initial covering $\{U_{i}\}$ so that each 
set $U_{i}$ is invariant under this generalized action. 

It is shown in~\cite{LT3} that these methods allow one to carry through the 
arguments in~\S\ref{ss:sket}.  Hence 
Proposition~\ref{prop:HZ1} holds for general 
$M$.  

Once we have this powerful method there is no need to cling to all the
special conditions that we put on $F$ that adapted it to the fibration
on $M\times S^2$.  For the argument to make sense, we need $F$ to be
constant and equal to its absolute maximum (resp. minimum) in a
neighborhood of one fiber and to assume its absolute minimum
(resp. maximum) at some point that plays the role of $p_{\infty}$.
The other important condition is that $F$ be slow.  Thus $F$ is
admissible in that it belongs to the set ${\cal H}_{ad}'(M\times
S^{2})$ defined in \S1. Using the methods of Liu--Tian to 
regularize the closure of the trajectory space $\Cc$ in $V =M\times S^2$
for these more general functions $F$, we obtain  
the following result.

\begin{prop}\label{prop:ac}  Given any closed symplectic manifold 
$(M,\om)$ and any 
quasicylinder $(Q = M\times D, \Om)$ the  
capacity $c_{HZ}'$ satisfies the area-capacity inequality
$$
c_{HZ}'(Q, \Om) \le {\rm area\,} (Q, \Om).
$$
\end{prop}

 Proposition~\ref{prop:cylcap} clearly follows.

\end{document}